\documentclass[11pt]{amsart}

\usepackage{etex}
\usepackage{amsmath, amssymb}
\usepackage{array}
\usepackage[frame,cmtip,arrow,matrix,line,graph,curve]{xy}
\usepackage{graphpap, color, paralist, pstricks}
\usepackage[mathscr]{eucal}
\usepackage[pdftex]{graphicx}
\usepackage[pdftex,colorlinks,backref=page,citecolor=blue]{hyperref}
\usepackage{tikz}
\usepackage{setspace}

\numberwithin{equation}{section}

\newtheorem{Thm}{Theorem}[section]
\newtheorem{Thm1}{Theorem}
\newtheorem{Thm2}[Thm1]{Theorem}
\newtheorem{Lem}[Thm]{Lemma}
\newtheorem{prop}[Thm]{Proposition}
\newtheorem{Cor}[Thm]{Corollary}

\theoremstyle{definition}
\newtheorem{Def}[Thm]{Definition}
\newtheorem{Nota}[Thm]{Notation}
\newtheorem{Rem}[Thm]{Remark}
\newtheorem{Ex}[Thm]{Example}
\newtheorem{Ass}[Thm]{Assumption}

\newcommand{\cA}{\mathcal{A} }
\newcommand{\cB}{\mathcal{B} }
\newcommand{\cC}{\mathcal{C} }
\newcommand{\cF}{\mathcal{F} }
\newcommand{\cG}{\mathcal{G} }
\newcommand{\cE}{\mathcal{E} }

\newcommand{\cN}{\mathcal{N} }
\newcommand{\cL}{\mathcal{L} }
\newcommand{\cO}{\mathcal{O} }

\newcommand{\cR}{\mathcal{R} }

\newcommand{\cU}{\mathcal{U} }

\newcommand{\cc}{\mathbb{C}}

\newcommand{\pp}{\mathbb{P}}

\newcommand{\zz}{\mathbb{Z}}
\newcommand{\rr}{\mathbb{R}}
\newcommand{\Pic}{\mathrm{Pic} }
\newcommand{\End}{\mathrm{End} }

\newcommand{\Hom}{\mathrm{Hom}}

\newcommand{\bC}{\mathbf{C} }
\newcommand{\bD}{\mathbf{D} }

\newcommand{\bH}{\mathbf{H} }

\newcommand{\bM}{\mathbf{M} }
\newcommand{\bS}{\mathbf{S} }

\newcommand{\bR}{\mathbf{R} }

\newcommand{\bK}{\mathbf{K} }

\newcommand{\bN}{\mathbf{N} }
\newcommand{\Gr}{\mathrm{Gr}}
\newcommand{\sgr}{\mathrm{gr}}

\newcommand{\id}{\textrm{id}}

\newcommand{\Aut}{\mathrm{Aut}}
\newcommand{\EEnd}{\mathscr E\!nd}
\newcommand{\bfHom}{\mathbf{Hom}}
\newcommand{\pHom}{\mathrm{ParHom}}
\newcommand{\spHom}{\mathrm{SParHom}}
\newcommand{\pHHom}{\mathscr P\!ar\mathscr H\!om}
\newcommand{\spHHom}{\mathscr {SP}\!ar\mathscr H\!om}
\newcommand{\HHom}{\mathscr H\!om}

\newcommand{\SL}{\mathrm{SL}}
\newcommand{\PGL}{\mathrm{PGL}}
\newcommand{\GL}{\mathrm{GL}}

\newcommand{\Supp}{\mathrm{Supp}}

\newcommand{\coker}{\mathrm{coker \,}}
\newcommand{\codim}{\mathrm{codim \,}}
\newcommand{\Stab}{\mathrm{Stab} \,}
\newcommand{\Spec}{\mathrm{Spec} \,}

\newcommand{\Sp}{\mathrm{Span}}

\newcommand{\git}{/\!\!/}

\newcommand{\tA}{\tilde{A} }
\newcommand{\tB}{\tilde{B} }

\newcommand{\tU}{\tilde{U} }

\newcommand{\tcN}{\tilde{\cN} }
\newcommand{\tcE}{\tilde{\cE} }
\newcommand{\tcF}{\tilde{\cF} }
\newcommand{\cD}{\mathcal{D} }
\newcommand{\tV}{\tilde{V} }
\newcommand{\tcD}{\tilde{\cD} }
\newcommand{\tD}{\tilde{D} }
\newcommand{\tY}{\tilde{Y} }
\newcommand{\tsig}{\tilde{\Sigma} }

\newcommand{\tdel}{\tilde{\Delta} }
\newcommand{\tM}{\widetilde{M} }

\newcommand{\tbN}{\tilde{\mathbf{N}} }
\newcommand{\tbR}{\tilde{\mathbf{R}} }

\newcommand{\tcL}{\tilde{\mathcal{L}} }

\newcommand{\e}{\epsilon }
\newcommand{\s}{\sigma }
\newcommand{\g}{\gamma }
\newcommand{\hs}{\hat{\sigma} }
\newcommand{\he}{\hat{\epsilon} }
\newcommand{\hg}{\hat{\gamma} }

\newcommand{\fS}{\mathfrak{S} }
\newcommand{\fm}{\mathfrak{m} }
\newcommand{\fp}{\mathfrak{p} }

\newcommand{\oUp}{\overline{\Upsilon}}
\newcommand{\oPsi}{\overline{\Psi}}
\newcommand{\opi}{\overline{\pi}}

\newcommand{\pdeg}{\mathrm{pardeg}\;}
\newcommand{\rk}{\mathrm{rank}\,}
\newcommand{\pa}{\mathrm{par}}

\newcommand{\fh}{\mathfrak{h}}

\hypersetup{%
pdftitle={A desingularization of the moduli space of rank 2 Higgs bundles over a curve},%
pdfauthor={Sang-Bum Yoo},%
pdfkeywords={desingularization, Higgs bundle, parabolic Higgs bundle, specialization of $M(2)$},%
citecolor=blue,%
linkcolor=blue,%
}

\begin{document}
\title[A desingularization of moduli of Higgs]
{A desingularization of the moduli space of rank 2 Higgs bundles over a curve}

\author{Sang-Bum Yoo}
\address{\parbox{\linewidth}{Department of Mathematics Education, Gongju National University of Education,\\ Gongju-si, Chungcheongnam-do, 32553, Republic of Korea}}
\email{sbyoo@gjue.ac.kr}

\begin{abstract}
Let $X$ be a smooth complex projective curve of genus $g\geq 3$. Let $\bM_2$ be the moduli space of semistable rank $2$ Higgs bundles with trivial determinant over $X$. We construct a desingularization $\bS$ of $\bM_2$ as a closed subvariety of a moduli space. We prove that $\bS$ is a nonsingular variety containing the stable locus of $\bM_2$ as an open dense subvariety. On the other hand, there is another desingularization $\bK$ of $\bM_2$ obtained from Kirwan's algorithm. We show that $\bS$ can be obtained after two blow-downs of $\bK$.
\end{abstract}


\maketitle

\section{Introduction}
Recently the birational geometry of various moduli spaces have become one of the main interests in algebraic geometry. It can be applied to construct new examples or to compute various invariants and so forth. On the moduli space of Higgs bundles over an algebraic curve, there have been a few developments in this direction. We can find some of them in \cite{BPG}, \cite{H} and \cite{Me}.

Throughout this paper, $X$ denotes a smooth complex projective curve of genus $g\ge 3$.

We first review a construction of a desingularization of the moduli space of vector bundles over an algebraic curve, which was done by C.S. Seshadri. Let $\bN$ be the moduli space of semistable rank $2$ vector bundles over $X$  with trivial determinant. It is known that $\bN$ has singularities along $\bN\setminus\bN^s$ where $\bN^s$ is the stable locus of $\bN$. In \cite{Se77} and \cite{Se82}, a nonsingular variety $\tbN$ was constructed together with a canonical morphism $\tbN\to\bN$ which is an isomorphism on $\bN^s$. $\tbN$ is indeed a closed subvariety of the moduli space of stable rank $4$, degree $0$ parabolic vector bundles over $X$ with respect to a point on $X$ with sufficiently small weights.

Seshadri's approach mentioned above can be applied to the moduli space of Higgs bundles over an algebraic curve. Let $\bM_2$ be the moduli space of semistable rank $2$ Higgs bundles over $X$  with trivial determinant. It is also known that $\bM_2$ has singularities along $\bM_2\setminus\bM_2^s$ where $\bM_2^s$ is the stable locus of $\bM_2$ (See \cite{KY} or \cite{Sim} for the details). We construct a desingularization $\bS$ of $\bM_2$. Let $\bM_{4,(a_1,a_2)}^{\pa}$ be the moduli space of semistable rank $4$, degree $0$ parabolic Higgs bundles over $X$ with multiplicities $(1,3)$ and weights $(a_1,a_2)$ with respect to a point on $X$. We choose sufficiently small weights $(a_1,a_2)$ so that $\bM_{4,(a_1,a_2)}^{\pa}=\bM_{4,(a_1,a_2)}^{\pa,s}$ where $\bM_{4,(a_1,a_2)}^{\pa,s}$ is the stable locus of $\bM_{4,(a_1,a_2)}^{\pa}$ (Proposition \ref{a choice of weights}). $\bS$ is defined as a closed subvariety of $\bM_{4,(a_1,a_2)}^{\pa}$. Precisely, $\bS$ is the Zariski closure of $\bS'$ in $\bM_{4,(a_1,a_2)}^{\pa,s}$, where $$\bS'=\{[(E,\phi,s)]\in\bM_{4,(a_{1},a_{2})}^{\pa,s}|\End((E,\phi))\cong M(2), \det E\cong\cO_{X}\text{ and }\phi\text{ is traceless}\}.$$
Here $M(2)$ is the $\cc$-algebra of $2\times 2$ matrices with entries in $\cc$. We have the following result.

\begin{Thm1}[Theorem \ref{smoothness}]\label{main1}
$\bS$ is nonsingular. Furthermore, there is a canonical morphism $\pi_{\bS}:\bS\to\bM_2$ which is an isomorphism on $\bM_2^s$.
\end{Thm1}

For the proof of the first statement, we first show that there exists a bijection between $\bM_2^s$ and $\bS'$ (Theorem \ref{set bijection}). Next, we define a moduli functor $par-Higgs_{4,(a_{1},a_{2})}^{sp}$ of $\bS$, which is the subfunctor of the moduli functor $par-Higgs_{4,(a_{1},a_{2})}$ of $\bM_{4,(a_1,a_2)}^{\pa}$ parametrizing all families of stable parabolic Higgs bundles equipped with specializations of $M(2)$ on the endomorphism algebras of the underlying Higgs bundles. Roughly speaking, the specialization of $M(2)$ is a limit structure of $\cc$-algebra structures of $M(2)$ (See \S\ref{specialization} for the details). Then we show that $par-Higgs_{4,(a_{1},a_{2})}^{sp}$ is formally smooth (Corollary \ref{cor:1}). To see this, we consider the morphism of functors $\Phi:par-Higgs_{4,(a_{1},a_{2})}^{sp}\to\fS_2$ which maps a family of stable parabolic Higgs bundles to the family of endomorphism algebras of the underlying Higgs bundles, where $\fS_2$ is the functor parametrizing all families of specializations of $M(2)$. We show that $\Phi$ is formally smooth (Proposition \ref{rel smoothness}). Finally we see that $\bS$ is a closed subvariety of $\bM_{4,(a_1,a_2)}^{\pa}$ and is a fine moduli scheme of $par-Higgs_{4,(a_{1},a_{2})}^{sp}$ (Proposition \ref{prop:5}). Then we conclude that $\bS$ is nonsingular.

The proof of the second statement is contained in that of Theorem \ref{smoothness}.

There is an interesting application of Theorem \ref{main1} (Theorem \ref{smoothness}). Let $\bK$ be Kirwan's desingularization of $\bM_2$ (See \S\ref{Kirwan's desingularization} in this paper, \cite{K2} or \cite{KY} for the details). Following O'Grady's argument in \cite{OGr} and \cite{OGr97}, a nonsingular variety $\bK_\e$ can be obtained after two blow-downs of $\bK$. But there is no known information on the moduli theoretic meaning of $\bK_\e$. The first author of \cite{KL} asked how we can provide the moduli theoretic meaning of $\bK_\e$.  We prove the following result.

\begin{Thm2}[Theorem \ref{main}]\label{main2}
$\bK_{\e}\cong\bS$.
\end{Thm2}

A similar work was already successfully completed on $\bN$ in \cite{KL}.

Here is an outline of this paper. In \S\ref{Preliminaries}, we recall basic facts about Higgs bundles, Parabolic Higgs bundles, Specialization of $M(r)$ and Upper semicontinuity for the complex of flat families of coherent sheaves. In \S\ref{Higgs bundles and Parabolic Higgs bundles}, we show that there exists a bijective correspondence between stable Higgs bundles of rank $r$ and stable parabolic Higgs bundles of rank $r^2$ with an endomorphism algebra of the underlying Higgs bundles as a specialization of $M(r)$. In \S\ref{Kirwan's desingularization}, we introduce the Kirwan's desingularization of $\bM_2$. In \S\ref{The construction of a desingularization}, we construct $\bS$ and then prove Theorem \ref{main1}. In \S\ref{A comparison between S and the Kirwan's desingularization}, we state Theorem \ref{main2} precisely. In \S\ref{A proof of Theorem 5.1}, we prove Theorem \ref{main2}.

\subsection*{Acknowledgements}
First of all I wish to thank my supervisor Young-Hoon Kiem for suggesting problem, and for his
help and encouragement. I am grateful to Conjeevaram Srirangachari
Seshadri for his inspiring papers \cite{Se77} and \cite{Se82}. I
also thank Insong Choe for useful comments.

\section{Preliminaries}\label{Preliminaries}
In this section, we recall some basic facts which are useful throughout this paper.

\subsection{Higgs bundles}
Let $r$ be an integer with $r\geq 2$. A \textit{Higgs bundle} of rank $r$ is a pair of a rank $r$
vector bundle $E$ with trivial determinant and a section $\phi$ of
$H^0(X,\EEnd_0 E\otimes K_X)$, where $K_X$ is the canonical bundle of $X$
and $\EEnd_0 E$ denotes the traceless part of $\EEnd E$. To construct
the moduli space of Higgs bundles, a stability condition has to be
imposed. C.T. Simpson introduced the following (See \cite{Sim}).

\begin{Def}
\begin{enumerate}
\item A Higgs bundle $(E,\phi)$ is \emph{(semi)stable} if for any nonzero proper subbundle $F$
satisfying $\phi (F)\subset F\otimes K_X$, we have
$$\mu(F):=\displaystyle\frac{\deg F}{\rk F}(\leq)<0.$$
\item A Higgs bundle $(E,\phi)$ is \emph{polystable}
if it is either stable or a direct sum of stable Higgs bundles.
\end{enumerate}
\end{Def}

A \textit{morphism of Higgs bundles} from $(E,\phi)$ to $(F,\psi)$
is a morphism of vector bundles $g:E\rightarrow F$ such that $\psi g=(g\otimes\text{id}_{K_X})\phi$, denoted by
$g:(E,\phi)\rightarrow(F,\psi)$. Throughout this paper, the space of morphisms of Higgs bundles from $(E,\phi)$ to $(F,\psi)$ is denoted by $\Hom((E,\phi),(F,\psi))$. In particular, $\Hom((E,\phi),(E,\phi))$ is denoted by $\End((E,\phi))$. We call $(F,\psi)$ a \textit{Higgs subbundle} of $(E,\phi)$ if there is an injective morphism of Higgs bundles $(F,\psi)\hookrightarrow(E,\phi)$.

Let $S$ be a scheme of finite type over $\cc$. Let $p_{1}:X\times S\to X$ and $p_{2}:X\times S\to S$ be the projections onto the first and the second components respectively. A pair $(\cE,\varphi)$ of a vector bundle $\cE$ on $X\times S$ and a morphism $\varphi:\cE\rightarrow\cE\otimes p_{1}^{*}K_{X}$ of vector bundles is \textit{a family of (semi)stable Higgs bundles} of rank $r$ on $X\times S$ if the restrictions $(\cE_{s},\varphi_{s})$ to the geometric fibers $X\times\{s\}$ are (semi)stable Higgs bundles of rank $r$ on $X$.

Let $(\cE,\varphi)$ be a family of (semi)stable Higgs bundles of rank $r$ on $X\times S$ and let $\{(U_i,z_i)\}$ be a pair of an open cover and local coordinate of $X$ such that $\cE|_{U_i}\cong\cO_{U_{i}\times S}^{\oplus r}$ and $\varphi|_{U_i}=A_{i}dz_{i}$ where $A_{i}$ is an $r\times r$ matrix with entries in $H^0(U_{i}\times S, \cO_{U_{i}\times S})$. We call $\{\varphi|_{U_i}\}$ (or $\{A_i\}$) \textit{a local Higgs field} of $(\cE,\varphi)$ with respect to $\{(U_i,z_i)\}$.

The set of isomorphism classes of polystable Higgs bundles of rank $r$ with trivial determinant admits a structure of irreducible normal quasi-projective variety of dimension $(r^2-1)(2g-2)$ (\cite{Hit, Nit, Sim}) and we denote it by $\bM_r$. The locus $\bM_r^s$ of stable Higgs bundles in $\bM_r$ is a smooth
open dense subvariety of $\bM_r$.

By \cite{Sim}, any Higgs bundle can be identified with a coherent
sheaf on the total space of the canonical bundle
$|K_{X}|=\mathbf{Spec}(S^{\bullet}K_{X}^{\vee})$. Let
$Z=\mathbf{Proj}(S^{\bullet}(K_{X}^{\vee}\oplus\mathcal{O}_{X}))$ and let
$D=Z-|K_{X}|$ be the divisor at infinity. The
projection $|K_{X}|\rightarrow X$ extends to a map
$\pi:Z\rightarrow X$. Choose a positive integer $k$ so that
$\mathcal{O}_{Z}(1):=\pi^{*}\mathcal{O}_{X}(k)\otimes_{\mathcal{O}_{Z}}\mathcal{O}_{Z}(D)$
is ample on $Z$. In particular,
$\mathcal{O}_{|K_{X}|}(1)=\pi^{*}\mathcal{O}_{X}(k)$. Thus for
any coherent sheaf $\mathcal{E}$ on $Z$ such that
$\Supp(\mathcal{E})\cap D=\emptyset$,
$\chi(\mathcal{E}(m))=\chi((\pi_{*}\mathcal{E})(km))$ by the
projection formula and the fact that $\pi|_{|K_{X}|}$ is affine. Furthermore, by virtue of \cite{BNR}, we can describe $\Supp(\cE)$ explicitly as a spectral curve.

\begin{Thm}[Proposition 3.6 in \cite{BNR}, Proposition 6.1 in \cite{HP} and Lemma 6.8 in \cite{Sim}]\label{BNR}
\begin{enumerate}
\item\label{identification between Higgs and coherent} A Higgs bundle $(E,\phi)$ on $X$ is the same thing as a
coherent sheaf $\mathcal{E}$ on $Z$ such that $\pi_{*}\cE=E$ and
$\Supp(\mathcal{E})\cap D=\emptyset$. This identification gives an
equivalence of categories.
\item The notions of semistability and stability for a Higgs bundle
$(E,\phi)$ on $X$ are the same as the corresponding notions for
the coherent sheaf
$\mathcal{E}$ on $Z$ associated to $(E,\phi)$ in (\ref{identification between Higgs and coherent}).
\item Any coherent sheaf $\mathcal{E}$ on $Z$ associated to
$(E,\phi)$ is of pure dimension $1$.
\item\label{spectral curve} For $\displaystyle s=(s_i)\in\bigoplus_{i=2}^{r}H^0(X,K_{X}^i)$, there is a bijective correspondence between isomorphism classes of torsion free sheaves $\cE$ on $X_s$ of rank $1$ and isomorphism classes of Higgs bundles $(E,\phi)$ of rank $r$ with coefficients $s_i$ of the characteristic polynomial of $\phi$ where the correspondence is given by (\ref{identification between Higgs and coherent}) and the spectral curve $X_s$ is defined by
$$x^{r}+s_{2}x^{r-2}+s_{3}x^{r-3}+\cdots+s_{r}=0$$
for $x$ the tautological section of $\pi^{*}K_X$.
\end{enumerate}
\end{Thm}

The following fact will be used later in \S\ref{Higgs bundles and Parabolic Higgs bundles}.

\begin{Lem}\label{torsion sheaves on Z}
Let $\cE$ be the coherent sheaf of pure dimension $1$ on $Z$ corresponding to
$(E,\phi)$ such that $\Supp(\cE)\cap D=\emptyset$. Then
$$\pi_{*}(\iota_{f,*}(\cE|_{f}))=\iota_{x_{0},*}(E|_{x_{0}})$$
where $x_0$ is a point of $X$, $\iota_{x_{0}}:\{x_0\}\hookrightarrow X$ is the inclusion of $\{x_0\}$ into $X$, $f=\pi^{-1}(x_0)\setminus D$ and $\iota_{f}:f\hookrightarrow|K_{X}|$ is the inclusion of $f$ into $|K_{X}|$.

Furthermore $\Supp(\iota_{f,*}(\cE|_{f}))$ is a zero-dimensional subscheme of length $\rk(E)$ in $|K_{X}|$.
\end{Lem}
\begin{proof}
We may assume that $X=\text{Spec}(A)$, $|K_{X}|=\text{Spec}(B)$ and $\pi:\text{Spec}(B)\to\text{Spec}(A)$.
Let $M$ be $A$-module and $B$-module, that is, $\tM$ is $\tA$-module and $\pi_{*}\tB$-module. Let $\fp$ be a prime ideal of $A$ corresponding to $x_0$. On $\text{Spec}(A)$
$$\iota_{x_{0},*}\tM^{A}|_{x_0}=\widetilde{M\otimes_{A}(A_{\fp}/\fp A_{\fp})}^{A}$$
and on $\text{Spec}(B)$
$$\iota_{f,*}\tM^{B}|_{f}=\widetilde{M\otimes_{B}(B\otimes_{A}(A_{\fp}/\fp A_{\fp}))}^{B}\cong\widetilde{M\otimes_{A}(A_{\fp}/\fp A_{\fp})}^{B}.$$
This completes the proof of the first statement.

From Theorem \ref{BNR}-(\ref{spectral curve}), $\Supp(\cE)$ is the spectral curve $X_s$ given by a polynomial of degree $\rk(E)$ and $\cE$ is a torsion free sheaf on $\Supp(\cE)$ of rank $1$. Since
$$-c_{2}(\iota_{f,*}(\cE|_{f}))=\chi(\iota_{f,*}(\cE|_{f})(m))=\chi(\iota_{x_{0},*}(E|_{x_{0}})(km))=c_{1}(\iota_{x_{0},*}(E|_{x_{0}}))=\rk(E),$$
we see that $\dim \Supp(\iota_{f,*}(\cE|_{f}))=0$ and $\text{length}(\Supp(\iota_{f,*}(\cE|_{f})))=\rk(E)$. This completes the proof of the second statement.
\end{proof}

\subsection{Extensions of Higgs bundles}
For Higgs bundles $(E,\phi)$ and $(E',\phi')$, let $\bfHom((E',\phi'),(E,\phi))$ be the complex
$$\xymatrix{\HHom(E',E)\ar[r]&\HHom(E',E)\otimes K_X}$$
given by $\psi\mapsto(\psi\otimes\text{id}_{K_X})\phi'-\phi\psi$. Then we have the following.

\begin{Thm}[(3.2) in \cite{HT}]\label{def sp of Higgs}
$\Hom((E',\phi'),(E,\phi))=\bH^0\bfHom((E',\phi'),(E,\phi))$. The space of extensions of $(E',\phi')$ by $(E,\phi)$ is $\bH^1\bfHom((E',\phi'),(E,\phi))$.
\end{Thm}

Here an \textit{extension} of $(E',\phi')$ by $(E,\phi)$ is a Higgs bundle $(E'',\phi'')$ which fits into the following commutative diagram
$$\xymatrix{0\ar[r]&E\ar[r]\ar[d]^{\phi}&E''\ar[r]\ar[d]^{\phi''}&E'\ar[r]\ar[d]^{\phi'}&0\\0\ar[r]&E\otimes K_X\ar[r]&E''\otimes K_X\ar[r]&E'\otimes K_X\ar[r]&0}$$
where each row is a short exact sequence.

\subsection{Parabolic Higgs bundles}
Let $r,n$ be integers with $r\geq 2$ and $n\geq 2$. Let $x_{0}$ be
a point of $X$. A \textit{parabolic bundle} of rank $n$ and degree $0$ on $X$, denoted $E_{*}$, with multiplicities
$(m_{1},m_{2},\cdots,m_{l})$ and weights
$(a_{1},a_{2},\cdots,a_{l})$ is a vector bundle $E$ of rank $n$ and degree $0$
on $X$ together with a filtration
$$E|_{x_0}=F_{1}(E)\supset F_{2}(E)\supset\cdots\supset F_{l}(E)\supset F_{l+1}(E)=\{0\},$$
weights
$$0\leq a_{1}<a_{2}<\cdots<a_{l}<1$$
and multiplicities
$$m_{i}=\dim_{\cc}(F_{i}(E)/F_{i+1}(E)).$$
Then $\pdeg(E_{*}):=\displaystyle\deg(E)+\sum_{i=1}^{l}m_{i}\cdot a_{i}$ is called the \textit{parabolic degree} of $E_{*}$. Denote $\displaystyle\pa\mu(E_{*}):=\frac{\pdeg(E_{*})}{\rk(E)}$.

Let $E_{*}$ and $E'_{*}$ be parabolic bundles on $X$ with weights
$$0\leq a_{1}(E)<a_{2}(E)<\cdots<a_{l}(E)<1$$
and
$$0\leq a_{1}(E')<a_{2}(E')<\cdots<a_{l}(E')<1$$
respectively. A bundle morphism $f : E \to E'$ is called \emph{(strongly) parabolic} if $f(F_{i}(E)) \subset F_{j+1}(E')$ whenever $a_{i}(E) (\ge) > a_{j}(E')$. The sheaves of parabolic morphisms and strongly parabolic morphisms are denoted by $\pHHom(E_{*}, E'_{*})$ and $\spHHom(E_{*}, E'_{*})$ respectively. The spaces of their global sections are denoted by $\pHom(E_{*}, E'_{*})$ and $\spHom(E_{*}, E'_{*})$.

If $E_*$ is a parabolic bundle on $X$ with a filtration
$$E|_{x_0}=F_{1}(E)\supset F_{2}(E)\supset\cdots\supset F_{l}(E)\supset F_{l+1}(E)=\{0\}$$
and weights
$(a_{1},a_{2},\cdots,a_{l})$ and $F$ is a subbundle of $E$, then a \textit{parabolic subbundle} $F_*$ of $E_*$ is given by
the filtration
$$F|_{x_0}=F_{1}(F)\supset F_{2}(F)\supset\cdots\supset F_{l}(F)\supset F_{l+1}(F)=\{0\}$$
such that $F_i(F)=F|_{x_0}\cap F_i(E)$, and the weights
$$(b_{1},b_{2},\cdots,b_{l})$$
such that $b_{j}=a_{i}$ for the largest index $i$ satisfying $F_j(F)=F|_{x_0}\cap F_i(E)$.

\begin{Def}
\emph{A parabolic Higgs bundle} $(E_{*},\phi)$ of rank $n$ and degree $0$ is a pair of a parabolic bundle $E_{*}$ of rank $n$ and degree $0$, and $\phi\in\pHom(E_{*},(E\otimes K_X)_{*})$ satisfying that $(E\otimes K_X)_*$
is a parabolic bundle of rank $n$ and degree $0$ such that $F_{i}(E\otimes K_X)=F_{i}(E)$
with the same weights as $E_{*}$.
\end{Def}

\begin{Def}
A parabolic Higgs bundle $(E_{*},\phi)$ is said to be \emph{(semi)stable} if for every proper parabolic
subbundle $F_{*}$ of $E_{*}$ satisfying $\phi(F_*)\subset(F\otimes K_X)_*$, we have $\pa\mu(F_{*})(\leq)<\pa\mu(E_{*})$.
\end{Def}

Let $S$ be a scheme of finite type over $\cc$. From now on, assume that $l=2$ and $(m_1,m_2)=(1,n-1)$. We write $(E_*,\phi)$ as $(E,\phi,s)$ where $F_2(E)$ is the hyperplane defined by $s\in\pp(E|_{x_0}^{\vee})$. If $F_*$ is a parabolic subbundle of $E_*$, then we write $F_*$ as $(F,s_F)$ where $s_F=\left\{\begin{matrix}s\text{ if }F|_{x_0}\not\subset F_2(E)\\0\text{ if }F|_{x_0}\subset F_2(E)\end{matrix}\right.$. A triple $(\cE,\varphi,\s)$ of a vector bundle $\cE$ on $X\times S$, a morphism $\varphi:\cE\rightarrow\cE\otimes p_{1}^{*}K_{X}$ of vector bundles and a section $\s\in\pp(\cE|_{\{x_0\}\times S}^{\vee})$ is \textit{a family of (semi)stable parabolic Higgs bundles} of rank $n$ and degree $0$ on $X\times S$ if $\cE$ is flat over $S$, $\cE|_{\{x_{0}\}\times S}$ has a filtration
$$\cE|_{\{x_{0}\}\times S}=F_{1}(\cE)\supset F_{2}(\cE)\supset F_{3}(\cE)=0$$
in which $F_2(\cE)$ is a vanishing locus of $\s$ and $\ker(\cE\to\cE|_{\{x_{0}\}\times S}/F_{2}(\cE))$ is flat over $S$, and if the restrictions $(\cE_{t},\varphi_{t},\s_t)$ to the geometric fibers $X\times\{t\}$ are (semi)stable parabolic Higgs bundles of rank $n$ and degree $0$ on $X$.

\begin{Def}
For each scheme $S$ of finite type over $\cc$, set
$$par-Higgs_{n,(a_{1},a_{2})}(S)=\{(\cE,\varphi,\s)|(\cE,\varphi,\s)\text{ is a family of semistable parabolic}$$
$$\text{Higgs bundles of rank }n\text{ and degree }0\text{ on }X\times S\text{ with the following properties}\}/\sim$$
such that
\begin{enumerate}
\item[(i)] for every geometric point $t\in S$, the restriction $(\cE_{t},\varphi_{t},\s_t)$ is a semistable parabolic Higgs bundle of rank $n$ and degree $0$ on $X$ with weights $0\le a_{1}<a_{2}<1$ and multiplicities $(1,n-1)$,

\item[(ii)] $(\cE,\varphi,\s)\sim(\cE',\varphi',\s')$ if and only if $(\cE,\varphi,\s)\cong(\cE',\varphi',\s')\otimes p_{2}^*\cL$ for some invertible sheaf $\cL$ on $S$.
\end{enumerate}
\end{Def}

Here $par-Higgs_{n,(a_{1},a_{2})}$ is a contravariant functor from the category of schemes of finite type over $\cc$ to the category of sets. $par-Higgs_{n,(a_{1},a_{2})}^{s}$ denotes the subfunctor of $par-Higgs_{n,(a_{1},a_{2})}$ consisting of all families of stable parabolic Higgs bundles.

\begin{Thm}[Theorem 4.6 and Corollary 4.7 in \cite{Yok1}, Theorem 5.2 in \cite{Yok2}]\label{coarse moduli}
There exists an irreducible normal quasi-projective variety $\bM_{n,(a_{1},a_{2})}^{\pa}$ over $\cc$ such that $\bM_{n,(a_{1},a_{2})}^{\pa}$ is a coarse moduli scheme of $par-Higgs_{n,(a_{1},a_{2})}$.
\end{Thm}

For parabolic Higgs bundles $(E,\phi,s)$ and $(E',\phi',s')$, let $\bfHom((E,\phi,s),(E',\phi',s'))$ be the complex
$$\xymatrix{\pHHom((E,s),(E',s'))\ar[r]&\spHHom((E,s),(E',s'))\otimes K_X}$$
given by $\psi\mapsto(\psi\otimes\text{id}_{K_X})\phi-\phi'\psi$. Then we have the following result.

\begin{Thm}[\cite{BR}]\label{long exact seq}
There is a long exact sequence
$$0\to\bH^0(\bfHom((E,\phi,s),(E',\phi',s')))\to H^0(\pHHom((E,s),(E',s')))\to H^0(\spHHom((E,s),(E',s'))\otimes K)$$
$$\to\bH^1(\bfHom((E,\phi,s),(E',\phi',s')))\to H^1(\pHHom((E,s),(E',s')))\to H^1(\spHHom((E,s),(E',s'))\otimes K)$$
$$\to\bH^2(\bfHom((E,\phi,s),(E',\phi',s')))\to0.$$
\end{Thm}

\subsection{The set of specialization of $M(r)$}\label{specialization}
This subsection will be useful in \S4 and \S5. Let $M(r)$ be the $\cc$-algebra of $r\times r$ matrices with entries in $\cc$. Fix a nonzero element $e_0\in\cc^{r^2}$. Let $\cA(r)$ be the set of elements in $\Hom(\cc^{r^2}\otimes\cc^{r^2},\cc^{r^2})$ which gives us an algebra structure on $\cc^{r^2}$ with the identity element $e_0$. There is a subset of $\cA(r)$ which consists of algebra structures on $\cc^{r^2}$, isomorphic to the matrix algebra $M(r)$. Let $\cA_r$ be the Zariski closure of this subset. An element in $\cA_r$ is called a \textit{specialization of} $M(r)$. Note that there exists a locally free sheaf $W$ of $\cO_{\cA_r}$-algebras on $\cA_r$ such that for each $z\in\cA_r$, $W|_{z}\otimes\cc$ is the specialization of $M(r)$ represented by $z$.

A \textit{family of specializations of} $M(r)$ \textit{parametrized by a noetherian}
$\cc$-\textit{scheme} $T$ is an $\mathcal{O}_{T}-$algebra $B$ such that for all
$t\in T$ there is a neighborhood $T_1$ of $t$ and a morphism
$f:T_{1}\rightarrow \cA_r$ satisfying $f^{*}(W)\cong
B|_{T_1}$.

Let $\fS_r$ be the functor from the category of schemes of finite type over $\cc$ to the category of sets
which to each $\cc$-scheme $T$ associates the set of isomorphism
classes of families of specializations of $M(r)$ parametrized by
$T$. It is known that $\fS_r$ is represented by $\cA_r$ and $\cA_2$ is smooth (See \cite[Chapter 5-I]{Se82}).

\section{Higgs bundles and Parabolic Higgs bundles}\label{Higgs bundles and Parabolic Higgs bundles}
In this section, we show that there exists a bijection set-theoretically between $\bM_r^s$ and a subset of $\bM_{r^2,(a_1,a_2)}^{\pa,s}$.

\begin{prop}\label{a choice of weights}
There exist $a_1$ and $a_2$ such that whenever $n\leq r^2$
\begin{enumerate}
\item[(i)] $par-Higgs_{n,(a_{1},a_{2})}(\Spec\cc)=par-Higgs^{s}_{n,(a_{1},a_{2})}(\Spec\cc)$,
\item[(ii)] for any $(E,\phi,s)\in par-Higgs_{n,(a_{1},a_{2})}(\Spec\cc)$, $(E,\phi)$ is a semistable Higgs
bundle,
\item[(iii)] For any $(E,\phi,s)\in par-Higgs_{n,(a_{1},a_{2})}(\Spec\cc)$ and any $\phi$-invariant subbundle $F$ of $E$, we have
$$\mu(F)<0\Rightarrow \pa\mu((F,s_F))<\pa\mu((E,s)).$$
\end{enumerate}
\end{prop}
\begin{proof}
Take $a_2$ such that $\displaystyle a_{2}<\frac{1}{r^{2}}$.
\begin{enumerate}
\item[(iii)] Let $(E,\phi,s)$ be an element of $par-Higgs_{n,(a_{1},a_{2})}(\Spec\cc)$. Assume that $\mu(F)<0$ for any $\phi$-invariant subbundle
$F$ of $E$.

If $F|_{x_0}\not\subset F_{2}(E)$, then
\begin{equation}\label{eq:1}
\pa\mu((E,s))-\pa\mu((F,s_F))=-\mu(F)+\frac{(a_{2}-a_{1})(n-\rk(F))}{n\cdot
\rk(F)}>0.
\end{equation}

If $F|_{x_0}\subset F_{2}(E)$, then
\begin{equation}\label{eq:2}
\pa\mu((E,s))-\pa\mu((F,s_F))=-\mu(F)+\frac{a_{1}-a_{2}}{n}.
\end{equation}
If $\displaystyle -\mu(F)+\frac{a_{1}-a_{2}}{n}\leq 0$, then
$$0<-\deg(F)n\leq(a_{2}-a_{1})\rk(F)<\frac{\rk(F)}{r^2}\leq\frac{\rk(F)}{n}<1,$$
which is a contradiction. Thus $\displaystyle-\mu(F)+\frac{a_{1}-a_{2}}{n}>0$.\\
\item[(i)] Let $(E,\phi,s)$ be an element of $par-Higgs_{n,(a_{1},a_{2})}(\Spec\cc)$ and let $F$ be a $\phi$-invariant subbundle of $E$. Then we have $\pa\mu((E,s))-\pa\mu((F,s_F))\geq 0$. By (iii), we have only to consider the case $\deg(F)\ge0$.

If $F|_{x_0}\not\subset F_{2}(E)$ and $\pa\mu((E,s))-\pa\mu((F,s_F))=0$ for some $\phi$-invariant proper subbundle $F$ of $E$,
then $\displaystyle 0=-\mu(F)+\frac{(a_{2}-a_{1})(n-\rk(F))}{n\cdot \rk(F)}$. Thus since $\displaystyle a_{2}<\frac{1}{r^{2}}\le\frac{1}{n}$,
$$0=-\deg(F)n+(a_{2}-a_{1})(n-\rk(F))<-\deg(F)n+\frac{n-\rk(F)}{n}.$$
Since $\deg(F)\geq 0$, we have $\displaystyle 0\leq\deg(F)n<\frac{n-\rk(F)}{n}<1$, which implies $\deg(F)=0$. It means that $(a_{2}-a_{1})(n-\rk(F))=0$, which is impossible.

If $F|_{x_0}\subset F_{2}(E)$ and
$\pa\mu((E,s))-\pa\mu((F,s_F))=0$ for some $\phi$-invariant proper subbundle $F$ of $E$, then $\displaystyle 0=-\mu(F)+\frac{a_{1}-a_{2}}{n}.$
So $\displaystyle -\deg(F)n=(a_{2}-a_{1})\rk(F)>0$. But $\deg(F)\geq 0$, which is a contradiction.\\
\item[(ii)] Let $(E,\phi,s)$ be an element of $par-Higgs_{n,(a_{1},a_{2})}(\Spec\cc)$ and let $F$ be a $\phi$-invariant subbundle of $E$. Then we have $\pa\mu((E,s))-\pa\mu((F,s_F))\geq 0$. Assume that $\mu(F)>0$.

If $F|_{x_0}\not\subset F_{2}(E)$ for some $\phi$-invariant proper subbundle $F$ of $E$, then
$$-\mu(F)+\frac{(a_{2}-a_{1})(n-\rk(F))}{n\cdot\rk(F)}\geq 0.$$
Since $\displaystyle a_{2}<\frac{1}{r^{2}}\le\frac{1}{n}$, $\displaystyle 0<\deg(F)n\leq(a_{2}-a_{1})(n-\rk(F))<\frac{n-\rk(F)}{n}<1$, which is a contradiction.

If $F|_{x_0}\subset F_{2}(E)$ for some $\phi$-invariant proper subbundle $F$ of $E$, then $\displaystyle-\mu(F)+\frac{a_{1}-a_{2}}{n}\geq 0$. So $\displaystyle0<\mu(F)\leq\frac{a_{1}-a_{2}}{n}<0$, which is a contradiction.
\end{enumerate}
\end{proof}

Under the choice of weights of Proposition \ref{a choice of weights}, $\bM_{n,(a_{1},a_{2})}^{\pa}=\bM_{n,(a_{1},a_{2})}^{\pa,s}$.

Let $(E,\phi)$ be a semistable Higgs bundle of rank $n$ on $X$. It is well known that $(E,\phi)$ has a Jordan-H\"{o}lder filtration
$$E=E_{0}\supset E_{1}\supset\cdots\supset E_{k+1}=0$$
such that for all $i$ $E_{i}$ is $\phi-$invariant and
$(E_{i}/E_{i+1},\overline{\phi})$ are stable Higgs bundles with
$\mu(E_{i}/E_{i+1})=\mu(E)$.

\begin{Lem}\label{PS of Higgs}
Every semistable Higgs bundle $(E,\phi)$ contains a unique nontrivial
maximal $\phi-$invariant subbundle $PS(E)$ of $E$ such that
$(PS(E),\phi|_{PS(E)})$ is a polystable Higgs bundle and
$\mu(PS(E))=\mu(E)$.
\end{Lem}
\begin{proof}
In a Jordan-H\"{o}lder filtration of $(E,\phi)$, there always exists a $\phi$-invariant proper subbundle $E_{1}$ such that $(E_{1},\phi|_{E_{1}})$ is stable and $\mu(E_{1})=\mu(E)$. By Zorn's Lemma, there is at least one nontrivial
maximal $\phi-$invariant subbundle $PS(E)$ of $E$ such that $(PS(E),\phi|_{PS(E)})$ is a polystable Higgs bundle and $\mu(PS(E))=\mu(E)$. Assume that there are two maximal polystable Higgs bundle $(PS_{1}(E),\phi|_{PS_{1}(E)})$ and $(PS_{2}(E),\phi|_{PS_{2}(E)})$. And assume that the uniqueness has been proved for all semistable pair $(E',\phi')$ with $\rk(E')<\rk(E)$. Both of $(PS_{1}(E),\phi|_{PS_{1}(E)})$ and $(PS_{2}(E),\phi|_{PS_{2}(E)})$ contain the stable Higgs bundle $(E_{1},\phi|_{E_{1}})$. Note that the induced Higgs bundle $(E/E_{1},\phi_{E/E_{1}})$ is semistable with $\mu(E/E_{1})=\mu(E)$. Then by induction hypothesis,
$$(PS(E/E_{1}),\phi_{E/E_{1}}|_{PS(E/E_{1})})\cong(PS_{1}(E)/E_{1},\phi|_{PS_{1}(E)/E_{1}})$$
$$\cong(PS_{2}(E)/E_{1},\phi|_{PS_{2}(E)/E_{1}}).$$
Therefore $(PS_{1}(E),\phi|_{PS_{1}(E)})\cong(PS_{2}(E),\phi|_{PS_{2}(E)})$.
\end{proof}

\begin{Rem}[\cite{HL}, Lemma 1.5.5]\label{PS of coh}
Let $P_{\cF}$ be the reduced Hilbert polynomial of a coherent sheaf $\cF$ on $Z$. Let $\cE$ be the coherent sheaf on $Z$ corresponding to $(E,\phi)$ such that $\Supp(\cE)\cap D=\emptyset$ mentioned in Theorem \ref{BNR}. By the same idea of the proof of Lemma \ref{PS of Higgs}, we can show that $\cE$ contains a unique nontrivial maximal
subsheaf $PS(\cE)$ of $\cE$ such that $PS(\cE)$ is a polystable sheaf
and $P_{PS(\cE)}=P_{\cE}$.
\end{Rem}

Then we have the following series of consequences.

\begin{prop}\label{surjective}
Let $(E,\phi)$ be a semistable Higgs bundle of rank $n$ on $X$. Then there exists a section $s$ such that $(E,\phi,s)\in par-Higgs^{s}_{n,(a_{1},a_{2})}(\Spec\cc)$ if
and only if for all stable Higgs bundles $(F,\psi)$ such that
$\mu(F)=\mu(E)$ and all integers $t>\rk(F)$ there is no injective
morphism of Higgs bundles $i:(F^{\oplus t},\psi^{\oplus t})\hookrightarrow
(E,\phi)$.
\end{prop}
\begin{proof}
We follow the idea of the proof of \cite[Proposition 7, Chapter 5]{Se82}. The proof of the \textit{if} part is nontrivial. Suppose that for all stable Higgs bundles $(F,\psi)$ such that
$\mu(F)=\mu(E)$, all integers $t>0$ and all injective morphism of
Higgs bundles $i:(F^{\oplus t},\psi^{\oplus t})\hookrightarrow (E,\phi)$, we
have
$t\leq \rk(F)$.\\
It suffices to show that for all $\phi-$invariant subbundle $F$ of
$E$ such that $(F,\phi|_{F})$ is a stable Higgs bundle and $\mu(F)=\mu(E)$,
there exists a hyperplane $F_{2}(E)$ of $E|_{x_0}$ satisfying
that $F|_{x_0}\not\subset F_{2}(E)$ and
$\phi$ is parabolic.

Assume that the previous statement is true. If $F$ is a $\phi-$invariant proper subbundle of $E$ such that
$\mu(F)<\mu(E)$, then
$$\pa\mu((E,s))-\pa\mu((F,s_F))=-\mu(F)+\frac{(a_{2}-a_{1})(n-\rk(F))}{n\cdot \rk(F)}>0.$$
If $F$ is a $\phi-$invariant proper subbundle of $E$ such that
$\mu(F)=\mu(E)$ and $(F,\phi|_{F})$ is a semistable Higgs bundle, then there
exists a hyperplane $F_{2}(E)$ of $E|_{x_0}$ such that
$F|_{x_0}\not\subset F_{2}(E)$ and $\phi$ is parabolic, because $(F,\phi|_{F})$ has a stable Higgs subbundle $(F',\phi|_{F'})$ such that $\mu(F')=\mu(E)$. In this case, we have also
$$\pa\mu((E,s))-\pa\mu((F,s_F))=\frac{(a_{2}-a_{1})(n-\rk(F))}{n\cdot \rk(F)}>0.$$

It remains to find the hyperplane $F_{2}(E)$ of $E|_{x_0}$ satisfying
that $F|_{x_0}\nsubseteq F_{2}(E)$ and
$\phi$ is parabolic, which is equivalent to find the $\phi$-invariant subsheaf $\iota_{x_{0},*}(F_{2}(E))$ of the skyscraper sheaf $\iota_{x_{0},*}(E|_{x_{0}})$ such that $\iota_{x_{0},*}(F|_{x_0})\nsubseteq\iota_{x_{0},*}(F_{2}(E))$ where $\iota_{x_{0},*}:\{x_{0}\}\hookrightarrow X$ is the inclusion and $F_{2}(E)$ is a hyperplane of $E|_{x_{0}}$. Note that $\iota_{x_{0},*}(E|_{x_{0}})$ is of pure dimension $0$.

Let $\cE$ be the coherent sheaf of pure dimension $1$ on $Z$ corresponding to
$(E,\phi)$ such that $\Supp(\cE)\cap D=\emptyset$ mentioned in Theorem \ref{BNR}. Write $(PS(E),\phi|_{PS(E)})=\oplus_{i=1}^{m}(F_{i}^{\oplus s_i},\psi_{i}^{\oplus s_i})$ where $(F_{i},\psi_{i})$ is a stable pair satisfying $\mu(F_i)=\mu(E)$ and $(F_{i},\psi_{i})\neq(F_{j},\psi_{j})$ for $i\neq j$. Then, by Lemma \ref{PS of Higgs} and Remark \ref{PS of coh}, every stable Higgs subbundle $(F,\phi|_{F})$ of $(E,\phi)$ such
that $\mu(F)=\mu(E)$ is isomorphic to $(F_{k},\psi_{k})$ for some
$k$ and that the corresponding stable subsheaf $\cF$ of $\cE$ on $Z$ such that
$P_{\cF}=P_{\cE}$ is
isomorphic to $\cF_{k}$. We may assume that $(F,\phi|_{F})=(F_{k},\psi_{k})$ and $\cF=\cF_{k}$. Let $f=\pi^{-1}(x_{0})\setminus D$. By Lemma \ref{torsion sheaves on Z}, $\pi_{*}(\iota_{f,*}(\cE|_{f}))=\iota_{x_{0},*}(E|_{x_{0}})$ where $\iota_{f}:f\hookrightarrow|K_{X}|$ is the inclusion of $f$ and $\iota_{x_{0}}:\{x_0\}\hookrightarrow X$ is the inclusion of $\{x_0\}$. It means that $\iota_{f,*}(\cE|_{f})$ is the coherent sheaf of pure dimension $0$ on $Z$ corresponding to $(\iota_{x_{0},*}(E|_{x_{0}}),\iota_{x_{0},*}(\phi|_{x_{0}}))$. Let
$$S_{k}:=\{(H,\cF)\in
\mathbb{P}(\cE|_{f}^{\vee})\times
\mathbb{P}(\cc^{s_{k}})|\cF|_{\Supp(\iota_{f,*}(\cE|_{f}))}\subset H\}.$$
Let $Z_{k}:=pr_{1}(S_{k})$. Then $\codim_{\mathbb{P}(\cE|_{f}^{\vee})}Z_{k}=\codim_{\mathbb{P}(\cE|_{f}^{\vee})}pr_{1}(pr_{2}^{-1}(\cF)\cap S_{k})=\rk(F_{k})\geq s_{k}$. So we can take a hyperplane $H$ in $\displaystyle\bigcap_{k=1}^{m}(\mathbb{P}(\cE|_{f}^{\vee})\setminus Z_{k})$.
Let $\widetilde{F_{2}(E)}=\pi_{*}\iota_{f,*}(H)$. Let $\lambda:\cE\to\cE\otimes\pi^{*}K_{X}$ be the morphism given by the tautological section $\lambda$ of $\pi^{*}K_{X}$. By \cite[Proposition 3.6]{BNR}, $\phi=\pi_{*}(\lambda):E\to E\otimes K_{X}$. Since $\iota_{f,*}(H)$ is a coherent sheaf of $\cO_{Z}$-module such that $\Supp(\iota_{f,*}(H))\cap D=\emptyset$, $\pi_{*}\iota_{f,*}(H)$ is $S^{\bullet}(K_{X}^{\vee})$-module and then it induces $\pi_{*}(\lambda):\pi_{*}\iota_{f,*}(H)\to\pi_{*}\iota_{f,*}(H)\otimes K_{X}$ (\cite[Lemma 2.13]{Sim}). So $$\phi(\widetilde{F_{2}(E)})=\pi_{*}(\lambda)(\pi_{*}\iota_{f,*}(H))\subset\pi_{*}\iota_{f,*}(H)\otimes K_{X}=\widetilde{F_{2}(E)}\otimes K_{X}.$$
Hence $\widetilde{F_{2}(E)}$ is a $\phi-$invariant subsheaf of $\iota_{x_{0},*}(E|_{x_{0}}$) such
that $\iota_{x_{0},*}(F|_{x_0})\nsubseteq\widetilde{F_{2}(E)}$, which completes
the proof.
\end{proof}

\begin{prop}\label{definition and well-definedness}
Let $(E,\phi,s)$ be an element of $par-Higgs^{s}_{n,(a_{1},a_{2})}(\Spec\cc)$ and let $M(r)$ be the $\cc$-algebra of $r\times r$ complex matrices. Then $\dim_{\cc}\End((E,\phi))\leq n$. Moreover when $n=r^{2}$, $\End((E,\phi))\cong M(r)$ of $\cc$-algebras if and only if there exists a stable
Higgs bundle $(F,\psi)$ such that
$$(E,\phi)\cong(F^{\oplus \rk(F)},\psi^{\oplus \rk(F)}).$$
\end{prop}
\begin{proof}
We follow the idea of the proof of \cite[Proposition 8, Chapter 5]{Se82}.
\end{proof}

\begin{prop}\label{injective}
Let $(E_{1},\phi_{1},s_1),(E_{2},\phi_{2},s_2)$ be elements of
$$par-Higgs^{s}_{n,(a_{1},a_{2})}(\Spec\cc)$$
such that
$$\dim\End((E_{1},\phi_{1}))=\dim\End((E_{2},\phi_{2}))=n.$$ Then $(E_{1},\phi_{1},s_1)\cong(E_{2},\phi_{2},s_2)$ if and only if
$(E_{1},\phi_{1})\cong(E_{2},\phi_{2})$.
\end{prop}
\begin{proof}
We follow the idea of the proof of \cite[Proposition 9, Chapter 5]{Se82}. The \textit{only if} part is obvious. For the proof of the \textit{if} part, suppose that $(E,\phi)$ be a semistable Higgs bundle of rank $n$ on $X$ which corresponds to a semistable coherent
sheaf $\cE$ such that $\Supp(\cE)\cap D=\emptyset$. Choose two $\phi$-invariant
hyperplanes $H_1$ and $H_2$ of $\cE|_{f}$ that make $(E,\phi)$ stable
parabolic Higgs bundles $(E_{1},\phi_{1},s_1)$ and $(E_{2},\phi_{2},s_2)$
respectively. The
group $\Aut_{\mathcal{O}_Z}(\cE)=\Aut_{\mathcal{O}_X}((E,\phi))$ acts
on
$$\bigcap_{i=1}^{m}(\mathbb{P}(\cE|_{f}^{\vee})\setminus Z_{i}).$$
Since every stable parabolic Higgs bundle is simple, the stabilizer of a point in
$\displaystyle\bigcap_{i=1}^{m}(\mathbb{P}(\cE|_{f}^{\vee})\setminus Z_{i})$ is isomorphic to
$\cc^{\times}$. Since $\dim\Aut_{\mathcal{O}_Z}(\cE)=n$, all orbits in
$\displaystyle\bigcap_{i=1}^{m}(\mathbb{P}(\cE|_{f}^{\vee})\setminus Z_{i})$ have dimension
$n-1$. Thus $\Aut_{\cO_Z}(\cE)$ acts on
$\displaystyle\bigcap_{i=1}^{m}(\mathbb{P}(\cE|_{f}^{\vee})\setminus Z_{i})$ transitively.
Therefore $H_{1}=g^{*}H_{2}$ for some $g\in\Aut_{\cO_Z}(\cE)$, that is,
$H_1$ and $H_2$ give an isomorphism of parabolic Higgs bundles $(E_{1},\phi_{1},s_1)$ and $(E_{2},\phi_{2},s_2)$.
\end{proof}

As a result, we have the following result.

\begin{Thm}\label{set bijection}
Assume that $\displaystyle a_{2}<\frac{1}{r^{2}}$ as Proposition \ref{a choice of weights}.
Let
$$\bS':=\{[(E,\phi,s)]\in\bM_{r^{2},(a_{1},a_{2})}^{\pa}=\bM_{r^{2},(a_{1},a_{2})}^{\pa,s}|\End((E,\phi))\cong M(r),\det E\cong\cO_{X}\text{ and }\phi\text{ is traceless}\}.$$
Then there exists a bijection
$$\bM_{r}^{s}\to\bS',(F,\psi)\mapsto(F^{\oplus r},\psi^{\oplus r},s_{can})$$
set-theoretically where $s_{can}$ is a canonical section so that $(F^{\oplus r},\psi^{\oplus r},s_{can})\in\bM_{r^{2},(a_{1},a_{2})}^{\pa,s}$ from Proposition \ref{surjective}.
\end{Thm}
\begin{proof}
The definition and well-definedness is obtained from Proposition \ref{surjective}, Proposition \ref{definition and well-definedness} and Proposition \ref{injective}. The surjectivity follows from Proposition \ref{definition and well-definedness}. The injectivity follows from Proposition \ref{injective}.
\end{proof}

\section{Kirwan's desingularization of $\bM_2$}\label{Kirwan's desingularization}

In this section, we briefly show that $\bM_2$ is desingularized by three
blow-ups by Kirwan's algorithm in the sense of \cite{K2}. For more details, see \cite[Section 4]{KY}.

In \cite[Theorem 3.8 and Theorem 4.10]{Sim}, C.T. Simpson showed
that a good quotient $\bR\git \SL(2)$ is $\bM_2$ where $\bR$ is an irreducible normal
quasi-projective variety which represents the moduli functor which
parameterizes triples $(V,\phi,\beta)$ where $(V,\phi)$ is a
semistable Higgs bundle with $\det V\cong \cO_X$, $\mathrm{tr}
\phi=0$ and $\beta$ is an isomorphism
\[
\beta:V|_x\to \cc^2
\]
and He also showed that every point in $\bR$ is semistable with respect to the action of $\SL(2)$, that the closed orbits in $\bR$ correspond to polystable Higgs bundles $(V,\phi)$, that is, stable or $(V,\phi)=(L,\psi)\oplus(L^{-1},-\psi)$ for $L\in\Pic^{0}(X)$ and $\psi\in H^{0}(K_{X})$ and that the set $\bR^{s}$ of stable points with respect to the action of $\SL(2)$ is exactly the locus of stable Higgs bundles.

We need to review a deformation theory for Higgs bundles discussed in \cite[Section 10]{Sim} and \cite[Section 4]{KY}. $A^{i}$ denotes the sheaf of smooth $i$-forms on $X$. For a polystable Higgs bundle $(V,\phi)$, we have the following complex
\begin{equation}\label{Dol complex}
\xymatrix{0\ar[r]&\EEnd_{0}V\otimes A^{0}\ar[r]^{\overline{\partial}+\phi}&\EEnd_{0}V\otimes A^{1}\ar[r]^{\overline{\partial}+\phi}&\EEnd_{0}V\otimes A^{2}\ar[r]&0}.
\end{equation}
This complex induces a long exact sequence
\begin{equation}\label{les from Dol complex}
\xymatrix{0\ar[r]&T^{0}\ar[r]&H^{0}(\EEnd_{0}V)\ar[r]^{[\phi,-]}&H^{0}(\EEnd_{0}V\otimes K_{X})\ar[r]&&\\
\ar[r]&T^{1}\ar[r]&H^{1}(\EEnd_{0}V)\ar[r]^{[\phi,-]}&H^{1}(\EEnd_{0}V\otimes K_{X})\ar[r]&T^{2}\ar[r]&0
}
\end{equation}
where $T^{i}$ is the $i$-th cohomology of (\ref{Dol complex}). Here the Zariski tangent space of $\bM_{2}$ at a polystable Higgs bundle $(V,\phi)$ is isomorphic to $T^{1}$. To investigate the singularities of $\bM_{2}$, we must recall the following theorem proved by C.T. Simpson.

\begin{Thm}[Theorem 10.4 and Theorem 10.5 of \cite{Sim}]\label{quadratic cone singularity}
Let $C$ be the quadratic cone in $T^{1}$ defined by the map $T^{1}\to T^{2}$ which sends an $\EEnd_{0}V$-valued $1$-form $\eta$ to $[\eta,\eta]$. Let $y=(V,\phi,\beta)\in\bR$ be a point with closed orbit and $\bar{y}\in\bM_{2}$ the point coming from $y$. Then the formal completion $(\bR,y)^{\wedge}$ is isomorphic to the formal completion $(C\times\fh^{\perp},0)^{\wedge}$ where $\fh^{\perp}$ is the perpendicular space to the image of the evaluation at $x$ composed with $\beta$, $T^{0}\to H^{0}(\EEnd_{0}V)\to sl(2)$. Furthermore, if $Y$ is the \'{e}tale slice at $y$ of the $\SL(2)$-orbit in $\bR$, then $(Y,y)^{\wedge}\cong(C,0)^{\wedge}$ and $(\bM_{2},\bar{y})^{\wedge}=(Y/\!/\Stab(y),\bar{y})^{\wedge}\cong(C/\!/\Stab(y),v)^{\wedge}$ where $\Stab(y)$ is the stabilizer of $y$ and $v$ is the cone point of $C$.
\end{Thm}

Let $\Omega_{\bR}$ (resp. $\Omega_{\bM_{2}}$) be the locus of $(L,0)\oplus (L,0)$ for $L\cong L^{-1}$ in $\bR\setminus\bR^s$ (resp. in $\bM_{2}\setminus\bM_{2}^s$). Both $\Omega_{\bR}$ and $\Omega_{\bM_{2}}$ are isomorphic to the $\zz_2$-fixed point set $\zz_2^{2g}$ in $J:=\Pic^{0}(X)$ by the involution $L\mapsto L^{-1}$. By (\ref{les from Dol complex}), we have isomorphisms
$$T^{0}\cong H^{0}(\EEnd_{0}V)\cong sl(2),$$
$$T^{1}\cong H^{0}(\EEnd{0}V\otimes K_{X})\oplus H^{1}(\EEnd_{0}V)\cong H^0(K_X)\otimes sl(2)\oplus H^1(\cO_X)\otimes sl(2)$$
and
$$T^{2}\cong H^{1}(\EEnd_{0}V\otimes K_{X})\cong H^{1}(K_{X})\otimes sl(2)\cong sl(2).$$
Then by Theorem \ref{quadratic cone singularity}, the normal cone of $\Omega_{\bM_{2}}$ in $\bM_{2}$ is a locally trivial fibration over $\Omega_{\bM_{2}}$ with fiber
$$\oUp^{-1}(0)/\!/ \SL(2)$$
where $\oUp:H^0(K_X)\otimes sl(2)\oplus H^1(\cO_X)\otimes sl(2)\to H^1(K_X)\otimes sl(2)\cong sl(2)$ is the quadratic map given by $(\alpha\otimes a,\beta\otimes b)\mapsto(\alpha\cup\beta)\otimes[a,b]$, the Lie bracket of $sl(2)$ coupled with the perfect pairing $H^0(K_X)\oplus H^1(\cO_X)\to H^1(K_X)$.

Let $\Sigma_{\bR}$ (resp. $\Sigma_{\bM_{2}}$) be the locus of $(L,\psi)\oplus (L^{-1},-\psi)$ for $(L,\psi)\ncong (L^{-1},-\psi)$ in $\bR\setminus\bR^s$ (resp. in $\bM_{2}\setminus\bM_{2}^s$). $\Sigma_{\bM_{2}}$ is isomorphic to
$$J\times_{\zz_2}H^0(K_X)-\zz_2^{2g}\cong T^*J/\zz_2-\zz_2^{2g}$$
where $\zz_2$ acts on $T^{*}J$ by $(L,\psi)\mapsto(L^{-1},-\psi)$. Since $\Sigma_{\bR}$ is a $\pp
\SL(2)/\cc^{\times}$-bundle over $T^*J/\zz_2-\zz_2^{2g}$, it is smooth. Since $L\ncong L^{-1}$, $H^{0}(\EEnd_{0}V)\cong H^{0}(\cO_{X})$ and then $[\phi,-]:H^{0}(\EEnd_{0}V)\to H^{0}(\EEnd_{0}V\otimes K_{X})$ is zero. Then by (\ref{les from Dol complex}), we have isomorphisms
$$T^{0}\cong\ker([\phi,-]:H^{0}(\EEnd_{0}V)\to H^{0}(\EEnd_{0}V\otimes K_{X}))\cong H^{0}(\cO_{X})\cong\cc,$$
$$T^{1}\cong\coker([\phi,-]:H^{0}(\EEnd_{0}V)\to H^{0}(\EEnd_{0}V\otimes K_{X}))\oplus\ker([\phi,-]:H^{1}(\EEnd_{0}V)\to H^{1}(\EEnd_{0}V\otimes K_{X}))$$
$$\cong H^{0}(\EEnd_{0}V\otimes K_{X})\oplus H^{1}(\EEnd_{0}V)\cong[H^{1}(\cO_{X})\oplus H^{0}(K_{X})]\oplus[H^0(L^{-2}K_X)\oplus H^0(L^{2}K_X)\oplus H^1(L^2)\oplus H^1(L^{-2})]$$
and
$$T^{2}\cong\coker([\phi,-]:H^{1}(\EEnd_{0}V)\to H^{1}(\EEnd_{0}V\otimes K_{X}))\cong H^{1}(\EEnd_{0}V\otimes K_{X})\cong H^{1}(K_{X}).$$
Note that $H^{1}(\cO_{X})\oplus H^{0}(K_{X})$ is the Zariski tangent space of $T^*J/\zz_2-\zz_2^{2g}$ at any point. Then by Theorem \ref{quadratic cone singularity}, the normal cone of $\Sigma_{\bM_{2}}$ in $\bM_{2}$ is a locally trivial fibration over $\Sigma_{\bM_{2}}$ with fiber
\[ \oPsi^{-1}(0)/\!/\cc^{\times},\]
where $\oPsi:H^0(L^{-2}K_X)\oplus H^0(L^{2}K_X)\oplus H^1(L^2)\oplus H^1(L^{-2})\to H^1(K_X)$ is the quadratic map given by $(a,b,c,d)\mapsto a\cup c+b\cup d$ for $(L,\psi)\oplus (L^{-1},-\psi)\in T^*J/\zz_2-\zz_2^{2g}$ and for perfect pairings $\cup:H^0(L^{-2}K_X)\oplus H^1(L^2)\to H^1(K_X)$ and $\cup:H^0(L^{2}K_X)\oplus H^1(L^{-2})\to H^1(K_X)$. Here the definition of $\oPsi$ comes from the diagonal entry in $\left[\left(\begin{matrix}\ast&b\\a&-\ast\end{matrix}\right),\left(\begin{matrix}\ast'&c\\d&-\ast'\end{matrix}\right)\right]$.

On the other hand, we have a stratification of $\bM_2$;
\[\bM_{2}=\bM_{2}^s\sqcup (T^*J/\zz_2-\zz_2^{2g})\sqcup
\zz_2^{2g}.\]

Since we have identical singularities and stratification as in
O'Grady's case in \cite{OGr}, we can adapt his arguments to construct the Kirwan's desingularization $\bK$
of $\bM_2$ and its blow-downs.

Let $\tbR$ be the variety obtained by blowing up $\bR$ first along $\Omega_{\bR}$ and then along the strict transform of $\Sigma_{\bR}$. Let $\tbR^{ss}$ (resp. $\tbR^s$) be the locus of semistable (resp. stable) points of $\tbR$. Then
we have \begin{enumerate}
\item[(a)] $\tbR^{ss}=\tbR^s$,
\item[(b)] $\tbR^s$ is smooth.\end{enumerate}
By blowing up $\tbR^s$ one more time along the locus $\tdel$ of points with stabilizers
larger than the center $\zz_2$ of $\SL(2)$, we obtain a variety
$\hat{\bR}$ with a smooth orbit space
$$\bK:=\hat{\bR}/\SL(2)$$
obtained by blowing up $\bM_2$ first along
$\zz_2^{2g}$, second along the strict transform of $T^*J/\zz_2$
and third along the nonsingular subvariety $\tdel/\!/\SL(2)$ contained in the strict transform of the exceptional divisor of the first blow-up. Let
$$\pi_{\bK}:\bK\to \bM_2$$
be the composition of these three blow-ups. $\bK$ is called the
\textit{Kirwan's desingularization} of $\bM_2$.

\section{The construction of a desingularization of $\bM_2$}\label{The construction of a desingularization}
In this section, we construct a desingularization of $\bM_2$. Assume that $r=2$ and fix $\displaystyle a_{2}<\frac{1}{4}$ as Proposition \ref{a choice of weights}. Then $\bM_{4,(a_{1},a_{2})}^{\pa}=\bM_{4,(a_{1},a_{2})}^{\pa,s}$. Let $par-Higgs_{4,(a_{1},a_{2})}^{sp}$ be
the subfunctor of
$par-Higgs_{4,(a_{1},a_{2})}$ consisting of
all families of stable parabolic Higgs bundles such that for
each $\cc$-scheme $T$ and each $(\cE,\varphi,\s)\in
par-Higgs_{4,(a_{1},a_{2})}(T)$, $p_{T
*}\cE\text{nd}((\cE,\varphi))\in\fS_{2}(T)$, where $p_{T}:X\times T\to T$ is the projection.

\subsection{Formally smoothness of $par-Higgs_{4,(a_{1},a_{2})}^{sp}$}
In this subsection, we prove that $par-Higgs_{4,(a_{1},a_{2})}^{sp}$ is formally smooth. Let $\cA$ be the category of local artinian commutative
$\cc$-algebras with residue field $\cc$. Let $\cF$ and $\cG$ be covariant functors from $\cA$ to the category of sets.

\begin{Def}
\begin{enumerate}
\item A morphism of functors $f:\cF\to\cG$ is called \emph{formally smooth} if given any surjective homomorphism $p:A'\to A$ in $\cA$ with the kernel $I$ such that $\fm_{A'}I=0$ and any elements $\alpha\in\cF(A)$ and $\beta\in\cG(A')$ such that
$$f_{A}(\alpha)=\cG(p)(\beta)\in\cG(A),$$
there exists an element $\gamma\in\cF(A')$ such that
$$f_{A'}(\gamma)=\beta\in\cG(A')\text{ and }\cF(p)(\gamma)=\alpha\in\cF(A).$$

\item A functor $\cF$ is called \emph{formally smooth} if $\cF(p):\cF(A')\to\cF(A)$ is surjective for any small surjection $p:A'\to A$ in $\cA$.
\end{enumerate}
\end{Def}

From the definition of $par-Higgs_{4,(a_{1},a_{2})}^{sp}$, we have the morphism of functors
$$\Phi:par-Higgs_{4,(a_{1},a_{2})}^{sp}\to\fS_2$$
given by $\Phi_{T}:par-Higgs_{4,(a_{1},a_{2})}^{sp}(T)\to\fS_2(T),\quad(\cE,\varphi,\s)\mapsto p_{T
*}\EEnd((\cE,\varphi))$ for each $\cc$-scheme $T$.

We claim that $\Phi$ is formally smooth, which implies that $par-Higgs_{4,(a_{1},a_{2})}^{sp}$ is formally smooth as a corollary.

\begin{Lem}\label{end lf}
Let $par-Higgs_{4,(a_{1},a_{2})}^{\textrm{end lf}}$ be the subfunctor of $par-Higgs_{4,(a_{1},a_{2})}$ consisting of all families of stable parabolic Higgs bundles such that for each $\cc$-scheme $T$ and each $(\cE,\varphi,\s)\in par-Higgs_{4,(a_{1},a_{2})}(T)$, $p_{T*}\cE\text{nd}((\cE,\varphi))$ is locally free, where $p_{T}:X\times T\to T$ is the projection. Then $par-Higgs_{4,(a_{1},a_{2})}^{\textrm{end lf}}$ is an open subfunctor of $par-Higgs_{4,(a_{1},a_{2})}$.
\end{Lem}
\begin{proof}
For any $(\cE,\varphi,\s)\in par-Higgs_{4,(a_{1},a_{2})}(T)$, let $U=\{t\in T\;|\;(p_{T
*}\cE\text{nd}((\cE,\varphi)))_{t}\text{ is locally free}\}$ which is an open subscheme of $T$. Then we can see that for any $f:T'\to T$, we have
$$f^{*}(\cE,\varphi,\s)\in par-Higgs_{4,(a_{1},a_{2})}^{\textrm{end lf}}(T')\Longleftrightarrow f\text{ factors through }U\subset T.$$
\end{proof}

\begin{Lem}\label{open subfunctor is fine}
The open subfunctor $par-Higgs_{4,(a_{1},a_{2})}^{\textrm{end lf}}$ of $par-Higgs_{4,(a_{1},a_{2})}$ defined in Lemma \ref{end lf} is represented by an open subscheme of $\bM_{4,(a_{1},a_{2})}^{\pa}$.
\end{Lem}
\begin{proof}
By Proposition \ref{fine}, $par-Higgs_{4,(a_{1},a_{2})}$ is represented by $\bM_{4,(a_{1},a_{2})}^{\pa}$. Let $(\cE^{\textrm{univ}},\varphi^{\textrm{univ}},\s^{\textrm{univ}})$ be the universal family. Consider $U^{\textrm{univ}}=\{t\in \bM_{4,(a_{1},a_{2})}^{\pa}\;|\;(p_{\bM_{4,(a_{1},a_{2})}^{\pa}
*}\cE\text{nd}((\cE^{\textrm{univ}},\varphi^{\textrm{univ}})))_{t}\text{ is locally free}\}$ which is an open subscheme of $\bM_{4,(a_{1},a_{2})}^{\pa}$. We can see that $par-Higgs_{4,(a_{1},a_{2})}^{\textrm{end lf}}$ is represented by $U^{\textrm{univ}}$.
\end{proof}

\begin{Lem}\label{open subfunctor is smooth}
The open subfunctor $par-Higgs_{4,(a_{1},a_{2})}^{\textrm{end lf}}$ of $par-Higgs_{4,(a_{1},a_{2})}$ defined in Lemma \ref{end lf} is formally smooth.
\end{Lem}
\begin{proof}
Since the open immersion $i:U^{\textrm{univ}}\hookrightarrow\bM_{4,(a_{1},a_{2})}^{\pa}$ is smooth, the morphism of functors $par-Higgs_{4,(a_{1},a_{2})}^{\textrm{end lf}}\to par-Higgs_{4,(a_{1},a_{2})}$ induced from $i$ is formally smooth by Lemma \ref{open subfunctor is fine} and the relative version of the infinitesimal lifting property (Exercise 4.7 of \cite{Hart10}).

Since we have chosen $a_1$ and $a_2$ such that $\bM_{4,(a_{1},a_{2})}^{\pa}=\bM_{4,(a_{1},a_{2})}^{\pa,s}$ in Proposition \ref{a choice of weights} and $\bM_{4,(a_{1},a_{2})}^{\pa,s}$ is smooth,
$$par-Higgs_{4,(a_{1},a_{2})}(p):par-Higgs_{4,(a_{1},a_{2})}(\Spec(A'))\to par-Higgs_{4,(a_{1},a_{2})}(\Spec(A))$$
is surjective for any surjective homomorphism $p:A'\rightarrow A$ in $\cA$ with the kernel $I$ such that $\fm_{A'}I=0$ by \cite[Proposition 2.2, Theorem 2.4, Theorem 5.2]{Yok2}. Hence we get the conclusion.
\end{proof}

\begin{prop}\label{rel smoothness}
The morphism of functors
$$\Phi:par-Higgs_{4,(a_{1},a_{2})}^{sp}\to\fS_2$$
is formally smooth.
\end{prop}
\begin{proof}
Let $(E,\phi,s)\in par-Higgs_{4,(a_{1},a_{2})}^{sp}(\Spec(A))$ and $B'\in\fS_2(\Spec(A'))$ such that $p_{A*}\EEnd((E,\phi))\cong\iota^*B'$ where $p:A'\rightarrow A$ is a surjective homomorphism in $\cA$ with the kernel $I$ such that $\fm_{A'}I=0$, $\iota:\Spec(A)\hookrightarrow\Spec(A')$ is the inclusion and $p_{A}:X\times\Spec(A)\to\Spec(A)$ is the projection. Note that $(E,\phi,s)\in par-Higgs_{4,(a_{1},a_{2})}^{\textrm{end lf}}(\Spec(A))$. It suffices to show that there exists a family of parabolic Higgs bundles $(E',\phi',s')\in par-Higgs_{4,(a_{1},a_{2})}^{sp}(\Spec(A'))$ such that $p_{A'*}\EEnd((E',\phi'))\cong B'\text{ and }\iota^*(E',\phi',s')\cong(E,\phi,s)$, where $p_{A'}:X\times\Spec(A')\to\Spec(A')$ is the projection.

We first claim that if there exists $(E',\phi',s')\in par-Higgs_{4,(a_{1},a_{2})}^{\textrm{end lf}}(\Spec(A'))$ such that $\iota^*(E',\phi',s')\cong(E,\phi,s)$, then $p_{A'*}\EEnd((E',\phi'))\cong B'$. Since $B'=f^{*}W$ for some morphism $f:\Spec(A')\to\cA_{2}$, $B'$ is locally free, that is, projective. Since we have the following commutative diagram
$$\xymatrix{B'\ar[r]&B'\otimes_{A'}A\ar[r]\ar[d]^{\cong}&0\\
p_{A'*}\EEnd((E',\phi'))\ar[r]&p_{A'*}\EEnd((E',\phi'))\otimes_{A'}A\ar[r]&0},$$
we get a surjective morphism $s:B'\to p_{A'*}\EEnd((E',\phi'))$. Since $p_{A'*}\EEnd((E',\phi'))$ is locally free by the hypothesis and the rank of $B'$ is equal to the rank of $p_{A'*}\EEnd((E',\phi'))$, $s:B'\to p_{A'*}\EEnd((E',\phi'))$ must be an isomorphism.

Next we show that there exists $(E',\phi',s')\in par-Higgs_{4,(a_{1},a_{2})}^{\textrm{end lf}}(\Spec(A'))$ such that $\iota^*(E',\phi',s')\cong(E,\phi,s)$. This is an immediate consequence of Lemma \ref{open subfunctor is smooth}.
\end{proof}

\begin{Cor}\label{cor:1}
$par-Higgs_{4,(a_{1},a_{2})}^{sp}$ is formally smooth.
\end{Cor}
\begin{proof}
For any surjective
homomorphism $p:A'\rightarrow A$ in $\cA$ with the kernel
$I$ such that $\fm_{A}I=0$, we have the following commutative diagram
$$\xymatrix{par-Higgs_{4,(a_{1},a_{2})}^{sp}(\Spec(A'))\ar[r]^{par-Higgs_{4,(a_{1},a_{2})}^{sp}(p)}\ar[d]_{\Phi(\Spec(A'))}&par-Higgs_{4,(a_{1},a_{2})}^{sp}(\Spec(A))\ar[d]^{\Phi(\Spec(A))}\\
\fS_2(\Spec(A'))\ar[r]^{\fS_2(p)}&\fS_2(\Spec(A))}.$$
Since $\cA_2$ is smooth, $\fS_2(p)$ is surjective. Combining this with Proposition \ref{rel smoothness}, $par-Higgs_{4,(a_{1},a_{2})}^{sp}(p)$ is surjective.
\end{proof}

\subsection{Fine moduli scheme of $par-Higgs_{4,(a_{1},a_{2})}^{sp}$}

In this subsection, we construct a fine moduli scheme $\bS$ of $par-Higgs_{4,(a_{1},a_{2})}^{sp}$ as a closed subvariety of $\bM_{4,(a_{1},a_{2})}^{\pa}$. Then we conclude that $\bS$ is indeed a desingularization of $\bM_2$.

\begin{prop}\label{fine}
$\bM_{4,(a_{1},a_{2})}^{\pa}$ is a fine moduli scheme of $par-Higgs_{4,(a_{1},a_{2})}$.
\end{prop}
\begin{proof}
\textsf{Step 1.}
We first claim that $\bM_{4,(a_{1},a_{2})}^{\pa}=\bM_{4,(a_{1},a_{2})}^{\pa,s}=\cR^{s}/\SL_{\nu(m)}=\cR^{s}/\PGL_{\nu(m)}$ for some quasi-projective variety $\cR$ with $\SL_{\nu(m)}$-action and $m\gg0$, where $\nu(m)=nm+n(1-g)$ and $\cR^{s}/\SL_{\nu(m)}$ is the moduli space of $1$-stable parabolic Higgs bundles of degree $0$ and rank $4$ on $X$ with multiplicities $(1,3)$ and weights $(a_{1},a_{2})$ constructed in \cite[\S2]{Yok1}. It suffices to show that every stable parabolic Higgs bundle $(E,\phi,s)$ is $1$-stable (See \cite[Definition 1.7]{Yok1}), that is, every nonzero subbundle $F$ of $E$ has $\deg(F)\le\rk(F)$.

Let $F_2(E)$ be the hyperplane given by $s$. Let $F$ be a subbundle of $E$. If $F|_{x_0}\not\subset F_{2}(E)$, then
$$\pa\mu((E,s))-\pa\mu((F,s_F))=-\mu(F)+\frac{(a_{2}-a_{1})(4-\rk(F))}{4\cdot
\rk(F)}>0,$$
which implies that
$$\deg(F)<\rk(F)\frac{(a_{2}-a_{1})(4-\rk(F))}{4\cdot
\rk(F)}<\rk(F)\frac{(4-\rk(F))}{4}<\rk(F).$$
If $F|_{x_0}\subset F_{2}(E)$, then
$$\pa\mu((E,s))-\pa\mu((F,s_F))=-\mu(F)+\frac{a_{1}-a_{2}}{4}>0,$$
which implies that $\displaystyle\deg(F)<\rk(F)\frac{a_{1}-a_{2}}{4}<\rk(F)$.

\textsf{Step 2.} Let $\cU$ be the universal family on $X\times\cR$, $\cU^{s}=\cU|_{X\times\cR^{s}}$ and let $\pi_{X}:X\times\cR^{s}\to X$ and $\pi_{\cR^{s}}:X\times\cR^{s}\to\cR^{s}$ be the projections onto the first and the second factor respectively. Note that $\GL_{\nu(m)}$-action on $\cR^{s}$ lifts to an action on $\cU^{s}$ and $\lambda(\id)$ acts on $\cU$ by scalar multiplication by $\lambda$ in fibers. If we are given a line bundle $L$ on $\cR^{s}$ such that $\GL_{\nu(m)}$-action on $\cR^{s}$ lifts to an action on $L$ and $\lambda(\id)$ acts on $L$ by scalar multiplication by $\lambda$, then $\cU^{s}\otimes\pi_{\cR^{s}}^{*}L^{-1}$ descends to the desired universal family on $X\times\bM_{4,(a_{1},a_{2})}^{\pa,s}$ by the descent lemma due to Kempf.

We have only to construct $L$ mentioned above. Define $\chi(k)=-4(2g-2)+m_{k}$. Let $(E,\phi,s_{E})$ be a parabolic Higgs bundle in $par-Higgs_{4,(a_{1},a_{2})}^{s}$ and let $(H_{k},\psi_{k},s_{H_{k}})$ be a parabolic line Higgs bundle with weights $a_{k}$ for $k=1,2$. Then
$$\dim\bH^0(\bfHom((H_{k},\psi_{k},s_{H_{k}}),(E,\phi,s_{E})))-\dim\bH^1(\bfHom((H_{k},\psi_{k},s_{H_{k}}),(E,\phi,s_{E})))$$
$$+\dim\bH^2(\bfHom((H_{k},\psi_{k},s_{H_{k}}),(E,\phi,s_{E})))=\chi(k)$$
by Theorem \ref{long exact seq}. Let
$$L(k)=\det\bR^{0}\pi_{R^s,*}\bfHom(\pi_{X}^{*}(H_{k},\psi_{k},s_{H_{k}}),\cU^{s})\otimes\det\bR^{1}\pi_{R^s,*}\bfHom(\pi_{X}^{*}(H_{k},\psi_{k},s_{H_{k}}),\cU^{s})^{-1}$$
$$\otimes\det\bR^{2}\pi_{R^s,*}\bfHom(\pi_{X}^{*}(H_{k},\psi_{k},s_{H_{k}}),\cU^{s}).$$
$\GL_{\nu(m)}$-action on $\cU$ induces a $\GL_{\nu(m)}$-action on $L(k)$ and $\lambda(\id)$ acts on $L(k)$ by scalar multiplication by $\lambda^{\chi(k)}$. Since $\chi(1)$ and $\chi(2)$ are consecutive odd numbers, they are relatively prime. Thus there exist $c_{1},c_{2}\in\zz$ such that $c_{1}\chi(1)+c_{2}\chi(2)=1$. Hence $L:=L(1)^{c_{1}}\otimes L(2)^{c_{2}}$ is the desired line bundle on $\cR^{s}$.
\end{proof}

The following is the construction of $\bS$.

\begin{prop}\label{prop:5}
Let $\bS$ be the subset of
$\bM_{4,(a_{1},a_{2})}^{\pa}$ consisting of
stable parabolic Higgs bundles $(E,\phi,s)$ such that
$\End((E,\phi))$
is a specialization of $M(2)$, $\det E\cong\cO_{X}$ and $\phi$ is traceless. Then
\begin{enumerate}
\item $\bS$ is a closed subvariety of $\bM_{4,(a_{1},a_{2})}^{\pa}$.

\item $\bS$ is a fine moduli scheme of
$par-Higgs_{4,(a_{1},a_{2})}^{sp}$.
\end{enumerate}
\end{prop}
\begin{proof}
\begin{enumerate}
\item The proof is identical as
that of \cite[Theorem 15-(i), Chapter 5]{Se82}. Let $\bS'$ be the subset of
$\bS$ consisting of
stable parabolic Higgs bundles $(E,\phi,s)$
such that $\End((E,\phi))\cong M(2)$ and $\det E\cong\cO_{X}$. It suffices to show that
$\mathbf{S}=\overline{\mathbf{S}'}$.

Let $(E,\phi,s)$ be an element of $\mathbf{S}$. Let $\cE$ be the coherent sheaf of pure dimension $1$ on $Z$ corresponding to $(E,\phi)$ such that $\Supp(\cE)\cap D=\emptyset$ mentioned in Theorem \ref{BNR}.

Since $\cA_2=\overline{\cA'_2}$, there exists a
$\mathcal{O}_{\Spec(\cc[[T]])}-$algebra $B$ as an element of
$\fS_{2}(\Spec(\cc[[T]]))$ such that
$$B\otimes_{\cc[[T]]}\cc\cong\End((E,\phi))$$
and
$$B\otimes_{\cc[[T]]}\overline{\cc((T))}\cong M(2,\overline{\cc((T))})$$
where $\cc((T))$ is the field of fractions of $\cc[[T]]$ and $\overline{\cc((T))}$ is the algebraic closure of $\cc((T))$.

If we are given a flat family
$(\cF,\psi,\s)$ as an element of
$$par-Higgs_{4,(a_{1},a_{2})}^{sp}(\Spec(\cc[[T]]))$$
such that $B\cong p_{0 *}\EEnd((\cF,\psi))$ and
$i^{*}(\cF,\psi,\s)\cong(E,\phi,s)$, where
$p_{0}:X\times\Spec(\cc[[T]])\rightarrow\Spec(\cc[[T]])$ is the
projection,
$i:\Spec(\cc)\hookrightarrow\Spec(\cc[[T]])$ is the inclusion, it gives us $\bS\subset\overline{\bS'}$ by \cite[Theorem of 3.1]{GM}. It is obvious that $\overline{\mathbf{S}'}\subset\mathbf{S}$.

It remains to prove the existence of $(\cF,\psi,\s)$. Let $D_{n}:=\cc[[T]]/(T^n)$ for all integer $n\geq 1$. This is an element of $\cA$. By Proposition \ref{rel smoothness}, there exists a sequence
$(\cF_{n},\psi_{n},\s_n)_{n\geq 1}$ such that for each integer $n\geq 1$,
\begin{enumerate}
\item[(a)] $(\cF_{n},\psi_{n},\s_n)$ is an element of
$par-Higgs_{4,(a_{1},a_{2})}^{sp}(\Spec(D_n))$,

\item[(b)] there exists an isomorphism $g_{n}:i_{n}^{*}(\cF_{n+1},\psi_{n+1},\s_{n+1})\rightarrow(\cF_{n},\psi_{n},\s_n)$
where $i_{n}:\Spec(D_n)\hookrightarrow\Spec(D_{n+1})$ is the
inclusion,

\item[(c)] there exists an isomorphism $f_{n}:p_{n *}\EEnd((\cF_{n},\psi_{n}))\rightarrow j_{n}^{*}B$
such that $f_{n+1}$ induces $f_n$ via $g_n$ where
$p_{n}:X\times\Spec(D_n)\rightarrow\Spec(D_n)$ is the projection
and $j_{n}:\Spec(D_n)\hookrightarrow\Spec(\cc[[T]])$ is the
inclusion.
\end{enumerate}

Then there exists a family of Higgs bundles $(\cF,\psi)$ uniquely on
$X\times\Spec(\cc[[T]])$ such that, for each integer $n\geq 1$, we
get an isomorphism $h_{n}:j_{n}^{*}(\cF,\psi)\rightarrow(\cF_{n},\psi_{n})$
such that $h_{n+1}$ induces $h_n$ via $g_n$. Precisely, let $U=\Spec(A)$ be an affine open subset of $X$. $(\cF_{n}|_{U\times\Spec(D_n)},\psi_{n})_{n\geq 1}$ on $X$
corresponds to $(\widetilde{\cF}_{n}|_{\pi^{-1}(U)\times\Spec(D_n)
})_{n\geq 1}$ on $Z$ in the sense of Theorem \ref{BNR}. This
sequence is equivalent to the datum of a sequence $(M_n)_{n\geq
1}$ such that for each integer $n\geq 1$,
\begin{enumerate}
\item[(a)] $M_n$ is a finitely generated $A\otimes_{\cc}D_{n}$-module,
\item[(b)] there exists an isomorphism $M_{n+1}\otimes_{A\otimes_{\cc}D_{n+1}}(A\otimes_{\cc}D_{n})\rightarrow M_{n}$.
\end{enumerate}
Put $\displaystyle M=\lim_{\stackrel{\longleftarrow}{n}}M_{n}$ which is a finitely generated $A\otimes_{\cc}\cc[[T]]$-module. Then we get
a coherent sheaf $\widetilde{\cF}$ on
$Z\times\Spec(\cc[[T]])$ which corresponds to a pair $(\cF,\psi)$ as desired.

Sections of $\widetilde{\cF}_{n}^{\vee}|_{f\times\Spec(D_n)}$ induced from $\s_n$ can be extended to a section $\s$ of
$\widetilde{\cF}^{\vee}|_{\Spec(f\times\cc[[T]])}$ such that $(\cF,\psi,\s)\in par-Higgs_{4,(a_{1},a_{2})}(\Spec(\cc[[T]]))$, and isomorphisms $f_{n}:p_{n *}\EEnd((\cF_{n},\psi_{n}))\rightarrow j_{n}^{*}B$ gives an isomorphism $f:p_{0 *}\EEnd((\cF,\psi))\rightarrow B$ by taking inverse limits.

\item By Proposition \ref{fine}, $\bM_{4,(a_1,a_2)}^{par}$ is fine. Then we can adapt the proof of \cite[Theorem 15-(ii), Chapter 5]{Se82}.
\end{enumerate}
\end{proof}

Then we get the first main result.

\begin{Thm}[Theorem \ref{main1}]\label{smoothness}
$\bS$ is nonsingular. Furthermore, there is a canonical morphism $\pi_{\bS}:\bS\to\bM_2$ which extends the isomorphism $\xymatrix{\bS'\ar[r]^{\cong}&\bM_2^s}$ given by $(E,\phi,s)\mapsto(F,\psi)$ where $(F,\psi)$ is the uniquely determined stable Higgs bundle such that $(E,\phi)\cong(F,\psi)\oplus(F,\psi)$.
\end{Thm}
\begin{proof}
The first statement is an immediate consequence from Corollary \ref{cor:1} and Proposition \ref{prop:5}.

Let's prove the second statement. Let $g:\bM_{2}\to\bM_{4}$ be a morphism given by $\sgr((F,\psi))\mapsto\sgr((F,\psi)\oplus(F,\psi))$, where $\sgr((F,\psi))$ and $\sgr((F,\psi)\oplus(F,\psi))$ are graded objects associated to Jordan-H\"{o}lder filtrations of $(F,\psi)$ and $(F,\psi)\oplus(F,\psi)$ respectively. It is obvious that $g$ is injective. Since $g$ induces an isomorphism between $\bM_{2}^{s}$ and $g(\bM_{2}^{s})$, $g:\bM_{2}\to g(\bM_{2})$ is birational. Indeed, $g$ is a homeomorphism onto $g(\bM_{2})$ because $\bM_{2}\setminus\bM_{2}^{s}=T^{*}J/\zz_{2}$ (See \S\ref{Kirwan's desingularization}) is separated (i.e. the diagonal morphism $\Delta:T^{*}J/\zz_{2}\to T^{*}J/\zz_{2}\times T^{*}J/\zz_{2}$ is a closed immersion) and $g(\bM_{2})\setminus g(\bM_{2}^{s})$ is homeomorphic to $\Delta(T^{*}J/\zz_{2})$.

Since $g(\bM_{2})$ and $p(\bS)$ contain the common dense subset $g(\bM_{2}^{s})$ from Theorem \ref{set bijection}, we have $g(\bM_{2})=p(\bS)$, where $p$ denotes the canonical forgetful morphism $\bM_{4,(a_1,a_2)}^{par}\to\bM_{4}$. Let $p'=p|_{\bS}$.

Since $\bS$ is nonsingular, the continuous function $g^{-1}\circ p':\xymatrix{\bS\ar[r]^{p'}&g(\bM_{2})\ar[r]^{g^{-1}}&\bM_{2}}$ is a morphism $\pi_{\bS}:\bS\to\bM_2$ by GAGA theorem and Riemann's extension theorem (\cite[4.10]{Mum76}). It is an immediate consequence from the construction that $\pi_{\bS}$ induces an isomorphism between $\pi_{\bS}^{-1}(\bM_{2}^{s})$ and $\bM_{2}^{s}$.
\end{proof}

\section{A comparison between $\bS$ and the Kirwan's desingularization of $\bM_2$}\label{A comparison between S and the Kirwan's desingularization}

In this section we solve a comparison problem, naturally raised from Theorem \ref{smoothness}, between $\bS$ and the Kirwan's desingularization of $\bM_2$.

\subsection{A comparison between $\bS$ and $\bK$}

By Proposition \ref{contraction}, $\bK$ have consecutive contractions
$$f:\xymatrix{\bK\ar[r]^{f_{\s}}&\bK_{\s}\ar[r]^{f_{\e}}&\bK_{\e}}$$
which are blow-downs. The following is the second main result which is a solution of this comparison problem.

\begin{Thm}[Theorem \ref{main2}]\label{main}
$\bK_{\e}\cong\bS$.
\end{Thm}

\section{A proof of Theorem \ref{main}}\label{A proof of Theorem 5.1}

In this section we first construct the morphism $\rho:\bK\to\bS$, show that $\rho$ factors through $\bK_{\e}$, that is,
$$\xymatrix{\bK\ar[rr]^{\rho}\ar[dr]_{f}&&\bS\\
&\bK_{\e}\ar[ur]_{\rho_{\e}}&}$$
and then we finally conclude that $\bK_{\e}\cong\bS$.

\subsection{A morphism from $\bK$ to $\bS$}\label{construction of morphism K to S}

We begin with the following proposition, which will be useful later.

\begin{prop}\label{key of construction of morphism}
Let $x_{0}\in X$.
\begin{enumerate}
\item Let $(E,\phi)\rightarrow T\times X$ be a family of
semistable Higgs bundles of rank $4$ and degree $0$ on $X$
parametrized by a complex manifold $T$. Assume the following:
\begin{enumerate}
\item for any $t\in T$, any line bundle $L$ of degree $0$ on
$X$ and any $\psi\in H^{0}(K_{X})$, $(L,\psi)\oplus(L,\psi)$ is not isomorphic to a Higgs subbundle of
$(E|_{t\times X},\phi|_{t\times X})$

\item there is an open dense subset $T'$ of $T$ such that
$$End((E|_{t\times X},\phi|_{t\times X}))\cong M(2)$$
for any $t\in T'$.
\end{enumerate}
Then we have a holomorphic map $\tau:T\rightarrow\mathbf{S}$.

\item Suppose a holomorphic map $\tau:T\rightarrow\mathbf{S}$ is
given. Suppose that $T$ is an open subset of a nonsingular
quasi-projective variety $W$ on which a reductive group $G$ acts
such that every point in $W$ is stable and the (smooth) geometric
quotient $q:W\rightarrow W/G$ exists. Furthermore, assume that
there is an open dense subset $W'$ of $W$ such that whenever
$t_{1},t_{2}\in T\cap W'$ are in the same orbit, we have
$\tau(t_1)=\tau(t_2)$. Then $\tau$ factors through the image $\bar{T}$ of
$T$ in $W/G$, that is, there is a continuous map $\overline{T}\to\bS$ such that the following diagram
$$\xymatrix{T\ar[rr]^{\tau}\ar[rd]&&\bS\\
&\overline{T}\ar[ru]&}$$
commutes.
\end{enumerate}
\end{prop}
\begin{proof}
\begin{enumerate}
\item Let $(E_{t},\phi_{t})=(E|_{t\times X},\phi|_{t\times X})$. By (1)-(a) and Proposition \ref{surjective}, there exists a section $0\ne s_{t}\in $ such that $(E_{t},\phi_{t},s_{t})$ is a stable parabolic Higgs bundle for each $t\in T$. Since $\bM_{4,(a_{1},a_{2})}^{\pa}$ is a coarse moduli scheme of $par-Higgs_{4,(a_{1},a_{2})}$ by Theorem \ref{coarse moduli}, we obtain a holomorphic map $\tau:T\to\bM_{4,(a_{1},a_{2})}^{\pa}$. By (1)-(b), $\tau(T')\subset\bS'$ and then $\tau$ is a map into $\bS$.

\item This is an immediate consequence of the \'{e}tale slice theorem. In particular, $\overline{T}$ is an open subset of $W/G$ in the usual complex topology.
\end{enumerate}
\end{proof}

Since both $\bK$ and $\bS$ contain Zariski open dense subsets which are isomorphic to $\bM_{2}^s$, there is a rational map
$$\rho':\bK\dashrightarrow\bS.$$

\begin{prop}\label{K to S}
There is a birational morphism $\rho:\bK\to\bS$ that extends $\rho'$
\end{prop}
\begin{proof}
The strategy of the proof is similar to that of \cite{KL}. By GAGA theorem and Riemann's extension theorem (\cite[4.10]{Mum76}), it
is enough to prove that $\rho'$ is extended to a continuous map with respect to the usual complex topology. Note that
$\bM_2=\bR\git \SL(2)$ for some irreducible normal quasi-projective
variety $\bR$ mentioned in \S\ref{Kirwan's desingularization}. Here $\bR$ is a
fine moduli scheme of a moduli functor  (See \S\ref{Kirwan's desingularization} or \cite{Sim} for more details). Luna's slice theorem
tells us that for each point $x\in\bM_2\setminus\bM_2^s$, there is
an analytic submanifold $W$ of $\bR$ such that $W/H$ is analytically equivalent to a neighborhood of
$x$ in $\bM_2$, where $H$ is the stabilizer of a point in the intersection of $W$ and the closed orbit
represented by $x$. Further, we can see that Kirwan's desingularization
$\widetilde{W}\git H$ of $W\git H$ is a neighborhood of the
preimage of $x$ in $\bK$. By Proposition \ref{key
of construction of morphism}, it is sufficient to
construct a nice family of semistable Higgs bundles of rank $4$
parametrized by $\widetilde{W}$ by successive applications of
elementary modifications of Higgs bundles, beginning with a
family of rank $2$ Higgs bundles parametrized by $W$, which comes
from the universal family over $X\times\bR$. Then we get a classifying morphism
$\widetilde{W}\to\bS$ because $\bS$ is the fine moduli scheme.
Since this is invariant under the $H$-action, we have an induced morphism $\widetilde{W}\git H\to\bS$. Therefore, $\rho'$ is
extended to a neighborhood of the preimage of $x$ in $\bK$. The detail of the construction of the morphism $\widetilde{W}\to\bS$ can be found in \S\ref{middle stratum}--\S\ref{D2 cap(D1 cup D3) stratum}.
\end{proof}

\subsubsection{Points over the middle stratum}\label{middle stratum}
We first extend $\rho'$ to points over the middle stratum of $\bM_2$. Let $l=[(L,\psi)\oplus(L^{-1},-\psi)]\in
T^*J/\zz_{2}-\zz_2^{2g}\subset\bM_2$ and let $\cN_{l}=H^0(K_X)\oplus H^1(\cO_X)\oplus\oPsi^{-1}(0)$. By Theorem \ref{quadratic cone singularity}, the universal object over $X\times\bR$ induces a holomorphic family $(\cF,\varphi_{\cF})$ of semistable Higgs bundles over $X$ parametrized by a neighborhood $U_{l}$ of $0$ in $\cN_{l}$. The restriction of $(\cF,\varphi_{\cF})$ to $X\times(U_{l}\cap(H^0(K_X)\oplus H^1(\cO_X)))$ is of the form $(\cL,\psi_{\cL})\oplus(\cL^{-1},-\psi_{\cL})$ where $\bR(\cc^{\times})$ is a quasi-projective variety parametrizing degree $0$ Higgs line bundles $(L,\psi)$ with an isomorphism $\beta:L|_{x}\to\cc$ (See \cite{Sim}) and $(\cL,\psi_{\cL})$ is a family of Higgs bundles of rank $1$ and degree $0$ that comes from an \'{e}tale map between $H^0(K_X)\oplus H^1(\cO_X)$ and the slice in $\bR(\cc^{\times})$.

Let $\pi_{l}:\tcN_{l}\to\cN_{l}$ be the blow-up of $\cN_{l}$ along $H^0(K_X)\oplus H^1(\cO_X)$. Let $\tU_{l}=\pi_{l}^{-1}(U_l)\cap\tcN_{l}^{ss}$ and $D_l$ be the exceptional locus in $\tU_{l}$. Let $(\tcF,\varphi_{\tcF})$ and $(\tcL,\psi_{\tcL})$ be the pull-backs of $(\cF,\varphi_{\cF})$ and $(\cL,\psi_{\cL})$ to $\tU_{l}$ and $D_l$ respectively. Then there are surjective morphisms
$$(\tcF|_{D_l},\varphi_{\tcF|_{D_l}})\to(\tcL,\psi_{\cL}),\quad\quad\quad(\tcF|_{D_l},\varphi_{\tcF|_{D_l}})\to(\tcL^{-1},-\psi_{\cL}).$$

Let $\tcF'$ and $\tcF''$ be the kernels of $\tcF\to\tcF|_{D_l}\to\tcL,\quad\tcF\to\tcF|_{D_l}\to\tcL^{-1}$ respectively. Then $\varphi_{\tcF'}:\tcF'\to\tcF'\otimes p_{1}^{*}K_X$ and $\varphi_{\tcF''}:\tcF''\to\tcF''\otimes p_{1}^{*}K_X$ can be determined uniquely such that
$$\xymatrix{0\ar[r]&\tcF'\ar[r]\ar[d]^{\varphi_{\tcF'}}&\tcF\ar[r]\ar[d]^{\varphi_{\tcF}}&\tcL\ar[r]\ar[d]^{\psi_{\cL}}&0\\
0\ar[r]&\tcF'\otimes p_{1}^{*}K_X\ar[r]&\tcF\otimes p_{1}^{*}K_X\ar[r]&\tcL\otimes p_{1}^{*}K_X\ar[r]&0}$$
and
$$\xymatrix{0\ar[r]&\tcF''\ar[r]\ar[d]^{\varphi_{\tcF''}}&\tcF\ar[r]\ar[d]^{\varphi_{\tcF}}&\tcL^{-1}\ar[r]\ar[d]^{-\psi_{\cL}}&0\\
0\ar[r]&\tcF''\otimes p_{1}^{*}K_X\ar[r]&\tcF\otimes p_{1}^{*}K_X\ar[r]&\tcL^{-1}\otimes p_{1}^{*}K_X\ar[r]&0}$$
commute. Consider
$(\cE,\varphi_{\cE})=(\tcF',\varphi_{\tcF'})\oplus(\tcF'',\varphi_{\tcF''})$
over $X\times\tU_l$. The following lemma tells us that $(\cE,\varphi_{\cE})$ satisfies the assumptions of Proposition \ref{key of construction of morphism}.

\begin{Lem}\label{useful lemma for the assumption}
$(\cE,\varphi_{\cE})$ is a family of semistable Higgs bundles of rank $4$, degree $0$ over $X$ parametrized by $\tU_l$ such that
\begin{enumerate}
\item For each $t\in\tU_l$, $L'\in\Pic^0(X)$ and $\psi'\in H^0(K_X)$, $(L',\psi')\oplus(L',\psi')$ is not isomorphic to any Higgs subbundle of $(\cE,\varphi_{\cE})|_{X\times t}$.
\item $(\cE,\varphi_{\cE})|_{X\times(\tU_{l}-D_l)}\cong((\tcF,\varphi_{\tcF})\oplus(\tcF,\varphi_{\tcF}))|_{X\times(\tU_{l}-D_l)}$ and there is an open dense subset of $\tU_l$ where $\End((\cE,\varphi_{\cE})|_{X\times t})$ is a specialization of $M(2)$.
\item With respect to the action of $\cc^{\times}$ on $\tcN_{l}-D_l$, if $t_{1},t_{2}\in\tU_{l}-D_l$ lie in the same orbit, then $(\cE,\varphi_{\cE})|_{X\times t_1}\cong(\cE,\varphi_{\cE})|_{X\times t_2}$.
\end{enumerate}
\end{Lem}
\begin{proof}
$\cE$ is locally free of rank $4$ because $D_l$ is a smooth divisor in $\tU_l$. Let $(a,b,c,d,e,f)\in\cN_{l}=H^0(K_X)\oplus H^1(\cO_X)\oplus\oPsi^{-1}(0)\subset H^0(K_X)\oplus H^1(\cO_X)\oplus[H^0(L^{-2}K_X)\oplus H^0(L^{2}K_X)]\oplus[H^1(L^2)\oplus H^1(L^{-2})]$. Since the weights of the $\cc^{\times}$-action are $0,0,-2,2,2,-2$ respectively, we can see that $(\cF,\varphi_{\cF})|_{X\times(a,b,c,d,e,f)}$ is stable if and only if $(d,e)\neq(0,0)$ and $(c,f)\neq(0,0)$.

For the proofs of (2) and (3), we can follow the proofs of \cite[Lemma 4.2-(2),(3)]{KL}. We have only to prove (1). For $t\in\tU_{l}-D_l$, it is obvious. Let $t$ be the point in $D_l$ represented by a line $C=(a,b,zc,zd,ze,zf)$ parametrized by $z\in\cc$ for $(a,b)\in H^0(K_X)\oplus H^1(\cO_X)$, $(0,0)\neq(d,e)\in H^0(L^{2}K_X)\oplus H^1(L^2)$ and $(0,0)\neq(c,f)\in H^0(L^{-2}K_X)\oplus H^1(L^{-2})$. Let $C_0=C\cap U_l$. Restricting $U_l$ sufficiently, we can take an open covering $\{V_i\}$ of $X$ such that $\cF|_{V_{i}\times C_0}$ are trivial. We fix a trivialization for each $i$. Let $(L_{ab},\psi_{ab})=(\cL,\psi)|_{X\times(a,b,0,0,0,0)}$. Since $(\cF,\varphi_{\cF})|_{X\times0}\cong(L_{ab},\psi_{ab})\oplus(L_{ab}^{-1},-\psi_{ab})$, the transition matrices are of the form
\begin{equation}\label{transition}
\left( \begin{matrix}\lambda_{ij}&z\alpha_{ij}\\z\beta_{ij}&\lambda_{ij}^{-1}\end{matrix}\right)
\end{equation}
where $\lambda_{ij}|_{z=0}$ is the transition for $L_{ab}$. $\varphi_{\cF}|_{V_{i}\times C_0}$ are of the form
\begin{equation}\label{local Higgs field}
\left( \begin{matrix}p_i&z q_i\\z r_i&-p_i\end{matrix}\right)
\end{equation}
which satisfies
\begin{equation}\label{compatibility of local Higgs field}
\left( \begin{matrix}p_j|_{V_{ij}\times C_0}&z q_j|_{V_{ij}\times C_0}\\z r_j|_{V_{ij}\times C_0}&-p_j|_{V_{ij}\times C_0}\end{matrix}\right)\left( \begin{matrix}\lambda_{ij}&z\alpha_{ij}\\z\beta_{ij}&\lambda_{ij}^{-1}\end{matrix}\right)=\left( \begin{matrix}\lambda_{ij}&z\alpha_{ij}\\z\beta_{ij}&\lambda_{ij}^{-1}\end{matrix}\right)\left( \begin{matrix}p_i|_{V_{ij}\times C_0}&z q_i|_{V_{ij}\times C_0}\\z r_i|_{V_{ij}\times C_0}&-p_i|_{V_{ij}\times C_0}\end{matrix}\right)
\end{equation}
where $p_{i}|_{z=0}=\psi_{ab}|_{V_{i}}$ and $V_{ij}=V_{i}\cap V_{j}$. The cocycle condition of (\ref{transition}) implies that $\{\lambda_{ij}\alpha_{ij}|_{z=0}\}$ and $\{\lambda_{ij}^{-1}\beta_{ij}|_{z=0}\}$ are cocycles in $H^{1}(L_{ab}^{2})$ and $H^{1}(L_{ab}^{-2})$ respectively. Since (\ref{local Higgs field}) are cocycles, $\{q_{i}|_{z=0}\}$ and $\{r_{i}|_{z=0}\}$ are cocycles in $H^{0}(L_{ab}^{2}K_{X})$ and $H^{0}(L_{ab}^{-2}K_{X})$ respectively. Since
$(\cF,\varphi_{\cF})|_{X\times(a,b,zc,zd,ze,zf)}$ is stable for
$z\neq 0$, $\{(\lambda_{ij}\alpha_{ij}|_{z=0},q_{i}|_{z=0})\}$ and $\{(\lambda_{ij}^{-1}\beta_{ij}|_{z=0},r_{i}|_{z=0})\}$ are all
nonzero.

Let $(\cF',\varphi_{\cF'})$ be the kernel of
$(\cF,\varphi_{\cF})|_{X\times
C_{0}}\to(\cF,\varphi_{\cF})|_{X\times
0}\cong(L_{ab},\psi_{ab})\oplus(L_{ab}^{-1},-\psi_{ab})\to(L_{ab},\psi_{ab})$ and let $(\cF'',\varphi_{\cF''})$ be the kernel of
$(\cF,\varphi_{\cF})|_{X\times
C_{0}}\to(\cF,\varphi_{\cF})|_{X\times
0}\cong(L_{ab},\psi_{ab})\oplus(L_{ab}^{-1},-\psi_{ab})\to(L_{ab}^{-1},-\psi_{ab})$.
Let $(F',\phi_{F'})=(\cF',\varphi_{\cF'})|_{X\times 0}$ and
$(F'',\phi_{F''})=(\cF'',\varphi_{\cF''})|_{X\times 0}$. Then
$$(\tcF',\varphi_{\tcF'})|_{X\times t}\cong(F',\phi_{F'})\text{ and }(\tcF'',\varphi_{\tcF''})|_{X\times t}\cong(F'',\phi_{F''}).$$
Note that any section of $\cF'$ over $V_{i}\times C_0$ is of the form $(z s_{1},s_{2})$. Since
\[\left( \begin{matrix}s_{1}\\s_{2}\end{matrix}\right)\leftrightarrow\left( \begin{matrix}z s_{1}\\s_{2}\end{matrix}\right)\mapsto\left( \begin{matrix}\lambda_{ij}&z\alpha_{ij}\\z\beta_{ij}&\lambda_{ij}^{-1}\end{matrix}\right)\left( \begin{matrix}z s_{1}\\s_{2}\end{matrix}\right)=\left( \begin{matrix}z(\lambda_{ij}s_{1}+\alpha_{ij}s_{2})\\z^{2}\beta_{ij}s_{1}+\lambda_{ij}^{-1}s_{2}\end{matrix}\right)\leftrightarrow\left( \begin{matrix}\lambda_{ij}s_{1}+\alpha_{ij}s_{2}\\z^{2}\beta_{ij}s_{1}+\lambda_{ij}^{-1}s_{2}\end{matrix}\right),\]
the transition for $\cF'$ is
\[\left( \begin{matrix}\lambda_{ij}&\alpha_{ij}\\z^{2}\beta_{ij}&\lambda_{ij}^{-1}\end{matrix}\right).\]
Since
\[\left( \begin{matrix}s_{1}\\s_{2}\end{matrix}\right)\leftrightarrow\left( \begin{matrix}z s_{1}\\s_{2}\end{matrix}\right)\mapsto\left( \begin{matrix}p_{i}&zq_{i}\\zr_{i}&-p_{i}\end{matrix}\right)\left( \begin{matrix}z s_{1}\\s_{2}\end{matrix}\right)=\left( \begin{matrix}z(p_{i}s_{1}+q_{i}s_{2})\\z^{2}r_{i}s_{1}-p_{i}s_{2}\end{matrix}\right)\leftrightarrow\left( \begin{matrix}p_{i}s_{1}+q_{i}s_{2}\\z^{2}r_{i}s_{1}-p_{i}s_{2}\end{matrix}\right),\]
$\varphi_{\cF'}|_{V_{i}\times C_0}$ is
\[\left( \begin{matrix}p_{i}&q_{i}\\z^{2}r_{i}&-p_{i}\end{matrix}\right).\]
By (\ref{compatibility of local Higgs field}),
\[\left( \begin{matrix}p_j|_{V_{ij}\times C_0}&q_j|_{V_{ij}\times C_0}\\z^{2}r_j|_{V_{ij}\times C_0}&-p_j|_{V_{ij}\times C_0}\end{matrix}\right)\left( \begin{matrix}\lambda_{ij}&\alpha_{ij}\\z^{2}\beta_{ij}&\lambda_{ij}^{-1}\end{matrix}\right)=\left( \begin{matrix}\lambda_{ij}&\alpha_{ij}\\z^{2}\beta_{ij}&\lambda_{ij}^{-1}\end{matrix}\right)\left( \begin{matrix}p_i|_{V_{ij}\times C_0}&q_i|_{V_{ij}\times C_0}\\z^{2}r_i|_{V_{ij}\times C_0}&-p_i|_{V_{ij}\times C_0}\end{matrix}\right).\]
Thus $(F',\phi_{F'})$ fits into a commutative diagram
$$\xymatrix{0\ar[r]&L_{ab}\ar[r]\ar[d]^{\psi_{ab}}&F'\ar[r]\ar[d]^{\phi_{F'}}&L_{ab}^{-1}\ar[r]\ar[d]^{-\psi_{ab}}&0\\
0\ar[r]&L_{ab}\otimes K_X\ar[r]&F'\otimes
K_X\ar[r]&L_{ab}^{-1}\otimes K_X\ar[r]&0}$$ where each row is
a short exact sequence and the extension class is $\{(\lambda_{ij}\alpha_{ij}|_{z=0},q_{i}|_{z=0})\}$, that is, $(F',\phi_{F'})$
is a nonsplit extension of $(L_{ab}^{-1},-\psi_{ab})$ by
$(L_{ab},\psi_{ab})$. Similarly $(F'',\phi_{F''})$ is also a
nonsplit extension of $(L_{ab},\psi_{ab})$ by
$(L_{ab}^{-1},-\psi_{ab})$ whose extension class is $\{(\lambda_{ij}^{-1}\beta_{ij}|_{z=0},r_{i}|_{z=0})\}$. Hence we conclude that $(\cE,\phi_{\cE})|_{X\times t}=(F',\phi_{F'})\oplus(F'',\phi_{F''})$ does
not contain a Higgs subbundle isomorphic to $(L',\psi')\oplus(L',\psi')$ for
any $L'\in\Pic^0(X)$ and $\psi'\in H^0(K_X)$.
\end{proof}
It follows from Proposition \ref{key of construction of morphism} and Luna's slice theorem that $\rho'$ extends continuously
to a neighborhood of the points in $\bK$ lying over $l$. Since there is at most one continuous extension, the extensions for various points
$l$ in the middle stratum $T^*J/\zz_{2}-\zz_2^{2g}$ are compatible. Thus $\rho'$ is extended to all the points in $\bK$ except those
over the deepest strata $\zz_2^{2g}$.

\subsubsection{Points over the deepest strata}\label{deepest stratum} We next extend $\rho'$ to the points over the deepest strata $\zz_2^{2g}$. We may consider only the
points in $\bK$ over $0=[(\cO,0)\oplus(\cO,0)]\in\zz_2^{2g}$. Let $\cN=\oUp^{-1}(0)$. Since it follows from Theorem \ref{quadratic cone singularity} that a
neighborhood of $[(\cO,0)\oplus(\cO,0)]$ in $\bM_2$ is
analytically equivalent to a neighborhood of the vertex
$\overline{0}$ in the cone $\cN/\!/\SL(2)$, a neighborhood of
the preimage of $[(\cO,0)\oplus(\cO,0)]$ in $\bK$ is bihomolorphic
to an open subset of the desingularization $\tcN/\!/\SL(2)$,
obtained by three blow-ups from $\cN/\!/\SL(2)$, which will be constructed below. Then it is enough to construct a holomorphic
map from a neighborhood $\tV$ of the preimage of $\overline{0}$ in
$\tcN/\!/\SL(2)$ to $\bS$.

Let $\Sigma=\SL(2)\{[H^0(K_X)\oplus H^1(\cO_X)]\otimes\left( \begin{matrix}1&0\\0&-1\end{matrix}\right)\}$. Let $\pi_{1}:\cN_{1}\to\cN$ be the blow-up of $\cN$ at $0$ with the exceptional divisor $\cD_1^{(1)}$. Let
$$\Delta=\SL(2)\pp\{\left( \begin{matrix}0&b\\c&0\end{matrix}\right)|b,c\in H^0(K_X)\oplus H^1(\cO_X),<b,c>=0\}$$
as a subset of $\cD_1^{(1)}$, where $<,>$ is the perfect pairing on $H^0(K_X)\oplus H^1(\cO_X)$. Let $\tsig$ be the strict transform of $\Sigma$ in $\cN_{1}$. Then the singular locus of $\cN_{1}/\!/\SL(2)$ is the quotient
$\Delta\cap\tsig/\!/\SL(2)$. Note that
\begin{equation}\label{intersection of D_1^(1) and tsig}\Delta\cap\tsig=\SL(2)\pp\{(H^0(K_X)\oplus H^1(\cO_X))\otimes\left( \begin{matrix}1&0\\0&-1\end{matrix}\right)\}=\cD_1^{(1)}\cap\tsig.\end{equation}

Let $\pi_{2}:\cN_{2}\to\cN_{1}$ be the blow-up of $\cN_{1}$ along $\tsig$ with the exceptional divisor $\cD_2^{(2)}$. Let
$\cD_2^{(1)}$ be the strict transform of $\cD_1^{(1)}$. The singular locus of $\cN_{2}/\!/\SL(2)$ is the quotient $\tdel/\!/\SL(2)$, where $\tdel$ is the strict transform of $\Delta$.

Let $\pi_{3}:\tcN=\cN_{3}\to\cN_{2}$ be the blow-up
of $\cN_2$ along $\tdel$ with the exceptional divisor $\tcD^{(3)}=\cD_3^{(3)}$ and let
$\tcD^{(1)}=\cD_3^{(1)},\tcD^{(2)}=\cD_3^{(2)}$ be the strict
transforms of $\cD_2^{(1)}$ and $\cD_2^{(2)}$ respectively. Let
$\pi:\tcN\to\cN$ be the composition $\pi_{1}\circ\pi_{2}\circ\pi_{3}$. Let $D_i^{(j)}$ be the quotient of $\cD_i^{(j)}$ in
$\cN_{i}/\!/\SL(2)$ for $1\leq i\leq 3$ and $1\leq j\leq i$.

It follows from a similar argument as in the case of the middle stratum that the universal family
over $X\times\bR$ induces a holomorphic family
$(\cF,\varphi_{\cF})$ of rank $2$ semistable Higgs bundles over
$X$ parametrized by an open neighborhood $U$ of $0$ in $\cN^{ss}$. Let
$V$ be the image of $U$ in $\cN/\!/\SL(2)$. Then $V$ is an open neighborhood of
$\overline{0}$ in $\cN/\!/\SL(2)$. Let
$U_{1}=\pi_{1}^{-1}(U)\cap\cN_{1}^{ss}$ and $V_1$ be the image of $U_1$ in $\cN_1/\!/\SL(2)$. It follows from
the commutative diagram
$$\xymatrix{\cN_1^{ss}\ar[r]\ar[d]^{\pi_1}&\cN_1/\!/\SL(2)\ar[d]^{\opi_1}\\
\cN^{ss}\ar[r]&\cN/\!/\SL(2)}$$
that we have $V_{1}=\opi_1^{-1}(V)$.

Let $U_2=\pi_{2}^{-1}(U_1)\cap\cN_{2}^{ss}$ and $V_2$ be the image
of $U_2$ in $\cN_2/\!/\SL(2)$. Then $V_{2}=\opi_2^{-1}(V_1)$ where
$\opi_{2}:\cN_2/\!/\SL(2)\to\cN_1/\!/\SL(2)$. Let
$\tU=\pi_{3}^{-1}(U_2)\cap\tcN^{ss}$ and $\tV$ be the image of
$\tU$ in $\tcN/\!/\SL(2)$. Then $\tV$ is smooth with
simple normal crossing divisors $\tD^{(1)},\tD^{(2)},\tD^{(3)}$
where $\tD^{(j)}=\tD_3^{(j)}$ and $\tD^{(2)}\cap\tV$ is denoted again by $\tD^{(2)}$.

$\rho'$ was already extended to the points over the middle
stratum, so we have a holomorphic map
$\rho':\tV-(\tD^{(1)}\cap\tD^{(3)})\to\bS$. Now $\rho'$ must be extended to $\rho:\tV\to\bS$.

\subsubsection{Points over $\tD^{(1)}-(\tD^{(2)}\cup\tD^{(3)})$}\label{D1-(D2 cup D3) stratum}
In this subsection, we will extend $\rho'$ to points in $\tV$ lying
over the quotient of $\cD_1^{(1)}-\Delta$.

We want to modify the pull-back of
$(\cF,\varphi_{\cF})\oplus(\cF,\varphi_{\cF})$ to
$U_{1}-\Delta\cup\tsig$ so that $\rho'$ is extended to a holomorphic
map near the quotient of $\cD_1^{(1)}-\Delta$ by Proposition
\ref{key of construction of morphism}.

Let $(\cF_{1},\varphi_{\cF_{1}})$ be the pull-back of
$(\cF,\varphi_{\cF})$ to $X\times U_{1}$ via $1_{X}\times\pi_{1}$. Since $(\cF,\varphi_{\cF})|_{X\times 0}$ is trivial, we have
$(\cF_{1},\varphi_{\cF_{1}})|_{X\times\cD_1^{(1)}}\cong(\cO,0)\oplus(\cO,0)$. Let $(\cF'_1,\varphi_{\cF'_1})$ (resp. $(\cF''_1,\varphi_{\cF''_1})$)
be the kernel of
$$(\cF_{1},\varphi_{\cF_{1}})\to(\cF_{1},\varphi_{\cF_{1}})|_{X\times\cD_1^{(1)}}\cong(\cO_{X\times\cD_1^{(1)}},0)\oplus(\cO_{X\times\cD_1^{(1)}},0)\to(\cO_{X\times\cD_1^{(1)}},0)$$
where the second arrow is the projection onto the first component (resp. the second component). A computation of transition matrices as in the proof of Lemma
\ref{useful lemma for the assumption} tells us that
$(\cF'_1,\varphi_{\cF'_1})|_{X\times t_1}$ and
$(\cF''_1,\varphi_{\cF''_1})|_{X\times t_1}$ are nonsplit
extensions of $(\cO,0)$ by $(\cO,0)$ for $t_{1}=\big(\left(
\begin{matrix}a&b\\c&-a\end{matrix}\right),\left(
\begin{matrix}d&e\\f&-d\end{matrix}\right)\big)\in\pp\cN=\cD_1^{(1)}$
with $(b,e)\neq(0,0)$ and $(c,f)\neq(0,0)$ in $\bH^{1}\bfHom((\cO,0),(\cO,0))=H^0(K_X)\oplus H^1(\cO_X)$.

Assume that $t_{1}\in\cD_1^{(1)}-\Delta$. Then
$(a,d),(b,e),(c,f)$ are linearly independent. In particular,
$(a,d),(b,e),(c,f)$ are all nonzero, which implies that
$(\cF'_1,\varphi_{\cF'_1})|_{X\times t_1}$ and
$(\cF''_1,\varphi_{\cF''_1})|_{X\times t_1}$ are nonsplit
extensions of $(\cO,0)$ by $(\cO,0)$ with extension classes $(b,e)$ and $(c,f)$ respectively.

The inclusion
$(\cF'_1,\varphi_{\cF'_1})\hookrightarrow(\cF_1,\varphi_{\cF_1})$
induces a homomorphism
$(\cF'_1,\varphi_{\cF'_1})|_{X\times\cD_1^{(1)}}\to(\cF_1,\varphi_{\cF_1})|_{X\times\cD_1^{(1)}}\cong(\cO,0)\oplus(\cO,0)$
whose image is the second factor $(\cO,0)$ and the kernel of this
homomorphism is $(\cO,0)$. Similarly,
$(\cO_{X\times\cD_1^{(1)}},0)$ is a Higgs subbundle of
$(\cF''_1,\varphi_{\cF''_1})$. Then we have a diagonal embedding of
$(\cO_{X\times\cD_1^{(1)}},0)$ into
$(\cF'_1,\varphi_{\cF'_1})\oplus(\cF''_1,\varphi_{\cF''_1})|_{X\times\cD_1^{(1)}}$.
Let $(\cE_1,\varphi_{\cE_{1}})$ be the kernel of
$$(\cF'_1\oplus\cF''_1,\varphi_{\cF'_1\oplus\cF''_1})\to(\cF'_1\oplus\cF''_1|_{X\times\cD_1^{(1)}},\varphi_{\cF'_1\oplus\cF''_1|_{X\times\cD_1^{(1)}}})\to(\cF'_1\oplus\cF''_1|_{X\times\cD_1^{(1)}}/\cO_{X\times\cD_1^{(1)}},\varphi_{(\cF'_1\oplus\cF''_1|_{X\times\cD_1^{(1)}}/\cO_{X\times\cD_1^{(1)}}}).$$

Adapting the argument of \cite[Section 4.3]{KL}, we can see that the transition for $\cE_{1}|_{X\times t_1}$ is
\begin{equation}\label{transition for cE_1|_{z=0}}
\left(
\begin{matrix}1&0&0&0\\
0&1&0&0\\
0&0&1&0\\
-2\alpha_{ij}|_{z=0}&-\beta_{ij}|_{z=0}&\gamma_{ij}|_{z=0}&1\end{matrix}\right).
\end{equation}
and $\varphi_{\cE_{1}}|_{X\times t_1}$ is locally
\begin{equation}\label{Higgs field of cE_1|_{z=0}}
\left(
\begin{matrix}0&0&0&0\\
0&0&0&0\\
0&0&0&0\\
-2p_i|_{z=0}&-q_i|_{z=0}&r_i|_{z=0}&0\end{matrix}\right).
\end{equation}
where the cocycles
$\{(\alpha_{ij}|_{z=0},p_{i}|_{z=0})\}$, $\{(\beta_{ij}|_{z=0},q_{i}|_{z=0})\}$,
$\{(\gamma_{ij}|_{z=0},r_{i}|_{z=0})\}$ represent the classes $(a,d),(b,e),(c,f)\in H^1(\cO_X)\oplus H^0(K_X)$ respectively.

It is easy to check that
$$\left(
\begin{matrix}0&0&0&0\\
0&0&0&0\\
0&0&0&0\\
-2p_j|_{z=0}&-q_j|_{z=0}&r_j|_{z=0}&0\end{matrix}\right)\left(
\begin{matrix}1&0&0&0\\
0&1&0&0\\
0&0&1&0\\
-2\alpha_{ij}|_{z=0}&-\beta_{ij}|_{z=0}&\gamma_{ij}|_{z=0}&1\end{matrix}\right)$$
$$=\left(\begin{matrix}1&0&0&0\\
0&1&0&0\\
0&0&1&0\\
-2\alpha_{ij}|_{z=0}&-\beta_{ij}|_{z=0}&\gamma_{ij}|_{z=0}&1\end{matrix}\right)\left(
\begin{matrix}0&0&0&0\\
0&0&0&0\\
0&0&0&0\\
-2p_i|_{z=0}&-q_i|_{z=0}&r_i|_{z=0}&0\end{matrix}\right).$$
Hence we have a filtration by Higgs subbundles
\begin{equation}\label{filtration of cE restricted to t1}
(\cE_1,\varphi_{\cE_1})|_{X\times t_1}=(E_4,\phi_{E_4})\supset(E_3,\phi_{E_3})\supset(E_2,\phi_{E_2})\supset(E_1,\phi_{E_1})\supset(E_0,\phi_{E_0})=0
\end{equation}
such that $(E_{i+1}/E_i,\phi_{E_{i+1}/E_i})\cong(\cO_X,0)$. Since $(c,f)\neq(0,0)$, the extension $(E_2,\phi_{E_2})$ of $(\cO,0)$ by $(E_1,\phi_{E_1})\cong(\cO,0)$ is nontrivial. Since $(b,e),(c,f)$ are linearly independent, an extension of $(\cO,0)$ by $(E_2,\phi_{E_2})$ is parametrized by $\bH^1\bfHom((\cO,0),(E_{2},\phi_{E_2}))$ which fits in the exact sequence
$$\xymatrix{\bH^0\bfHom((\cO,0),(\cO,0))\ar[r]^{(c,f)}&\bH^1\bfHom((\cO,0),(\cO,0))}$$
$$\xymatrix{\ar[r]&\bH^1\bfHom((\cO,0),(E_{2},\phi_{E_2}))\ar[r]&\bH^1\bfHom((\cO,0),(\cO,0))}$$
and $(E_3,\phi_{E_3})$ is the image of
$$(b,e)\in\bH^1(\bfHom((\cO,0),(\cO,0)))\cong H^1(\cO_X)\oplus H^0(K_X)$$
which is nonzero. Thus $(E_3,\phi_{E_3})$ is a nonsplit extension. Similarly $(E_4,\phi_{E_4})$ is a nonsplit extension since $(a,d),(b,e),(c,f)$ are linearly independent. Thus the condition $(1)-(a)$ of Proposition \ref{key of construction of morphism} holds for points in $\tU$ over $\cD_1^{(1)}-\Delta$. The other conditions of Proposition \ref{key of construction of morphism} obviously hold. Therefore $\rho'$ is extended to the points over the quotient of the points over $\cD_1^{(1)}-\Delta$.

\subsubsection{Points over $\tD^{(3)}-\tD^{(2)}$}\label{D3-D2 stratum} Let $t_{1}=\big(\left(
\begin{matrix}a&b\\c&-a\end{matrix}\right),\left(
\begin{matrix}d&e\\f&-d\end{matrix}\right)\big)\in\Delta-\tsig$. Then $(a,d),(b,e),(c,f)$ span 2-dimensional subspace of $H^1(\cO_X)\oplus H^0(K_X)$. The Higgs bundle $(\cE_{1},\varphi_{\cE_{1}})|_{X\times t_1}$ has transition matrices and local Higgs fields of the form (\ref{transition for cE_1|_{z=0}}) and (\ref{Higgs field of cE_1|_{z=0}}). A Higgs bundle of rank $1$ over $\Delta-\tsig$ is obtained from the family of linear relations of $(a,d),(b,e),(c,f)$. $(\cL_{1},0)$ denotes the pull-back of this Higgs bundle of rank $1$ to $X\times(\Delta-\tsig)$. Then there is an embedding of $(\cL_{1},0)\hookrightarrow(\cE_{1},\varphi_{\cE_{1}})|_{X\times(\Delta-\tsig)}$. Let
$(\cE_{3},\varphi_{\cE_{3}})$ (resp. $(\cL_{3},0)$) be the
pull-back of $(\cE_{1},\varphi_{\cE_{1}})$ (resp. $(\cL_{1},0)$)
to $\tU=U_{3}$ (resp. $\tcD^{(3)}-\tcD^{(2)}$).

Let $(\tcE,\varphi_{\tcE})$ be the kernel of
$$(\cE_{3},\varphi_{\cE_{3}})\to(\cE_{3}|_{X\times(\tcD^{(3)}-\tcD^{(2)})},\varphi_{\cE_{3}|_{X\times(\tcD^{(3)}-\tcD^{(2)})}})\to(\cE_{3}|_{X\times(\tcD^{(3)}-\tcD^{(2)})}/\cL_{3},\varphi_{\cE_{3}|_{X\times(\tcD^{(3)}-\tcD^{(2)})}/\cL_{3}}).$$

Let $t_{1}=\big(\left(
\begin{matrix}0&b\\c&0\end{matrix}\right),\left(
\begin{matrix}0&e\\f&0\end{matrix}\right)\big)\in\Delta-\tsig$
with $(b,e),(c,f)$ linearly independent for simplicity. (The general case is
obtained by conjugation.) Let $t_{3}\in\tcD^{(3)}-\tcD^{(2)}$ be a
(semi)stable point lying over $t_{1}$. A point $t_{3}\in\tcD^{(3)}$ represents a normal direction to
$\Delta$ at $t_{1}$. Choose a local parameter $z$ to the direction
such that $z=0$ represents $t_{1}$. We make local computations as in the previous subsection \S\ref{D1-(D2 cup D3) stratum}.

If $t_{3}$ represents a normal direction of $\Delta$ tangent to
$\tcD^{(1)}$, then the transition and local Higgs field of $(\tcE,\varphi_{\tcE})|_{X\times t_{3}}$ are of the form
\begin{equation}\label{transition for tcE 1}
\left(
\begin{matrix}1&0&0&0\\
0&1&0&0\\
0&0&1&0\\
-2\delta_{ij}|_{z=0}&-\beta_{ij}|_{z=0}&\gamma_{ij}|_{z=0}&1\end{matrix}\right)
\end{equation}
and
\begin{equation}\label{Higgs field of tcE 1}
\left(
\begin{matrix}0&0&0&0\\
0&0&0&0\\
0&0&0&0\\
-2s_i|_{z=0}&-q_i|_{z=0}&r_i|_{z=0}&0\end{matrix}\right)
\end{equation}
where $\{(\delta_{ij},s_{i})\}$ are some cocycles which give us a nonzero class $(g,h)\in H^0(K_X)\oplus H^1(\cO_X)$ at $z=0$ such that $(g,h),(b,e),(c,f)$ are linearly independent. The condition $(1)$ of Proposition \ref{key of construction of morphism} holds because the Higgs bundle has a filtration by Higgs subbundles as in (\ref{filtration of cE restricted to t1}).

If $t_{3}$ represents the direction normal to $\tcD^{(1)}$, then the transition and local Higgs field for
$(\tcE,\varphi_{\tcE})|_{X\times t_{3}}$ are of the form
\begin{equation}\label{transition for tcE 2}
\left(
\begin{matrix}1&0&0&0\\
\gamma_{ij}|_{z=0}&1&0&0\\
\beta_{ij}|_{z=0}&0&1&0\\
-2\delta_{ij}|_{z=0}&-\beta_{ij}|_{z=0}&\gamma_{ij}|_{z=0}&1\end{matrix}\right)
\end{equation}
and
\begin{equation}\label{Higgs field of tcE 2}
\left(
\begin{matrix}0&0&0&0\\
r_i|_{z=0}&0&0&0\\
q_i|_{z=0}&0&0&0\\
-2s_i|_{z=0}&-q_i|_{z=0}&r_i|_{z=0}&0\end{matrix}\right).
\end{equation}
Since the Higgs bundle has a filtration by Higgs subbundles as
in (\ref{filtration of cE restricted to t1}), $(\tcE,\varphi_{\tcE})|_{X\times
t_{3}}$ satisfies the condition (1) of Proposition \ref{key of
construction of morphism}.

Since the other conditions of Proposition \ref{key of construction of morphism} automatically hold on the stable part of $U$, $\rho'$ extends to the quotient of $\tU-\tcD^{(2)}$.

\subsubsection{Points over $\tD^{(2)}\cap(\tD^{(1)}\cup\tD^{(3)})$}\label{D2 cap(D1 cup D3) stratum}
In this subsection $\rho'$ will be finally extended to $\bK$ and then we complete the proof of Proposition \ref{K to S}. The slice theorem gives us a map $\tV\to\bK$ which is biholomorphic onto a neighborhood of the preimage of $[(\cO,0)\oplus(\cO,0)]$. So it is sufficient to construct a holomorphic map $\tV\to\bS$.

There is a commutative diagram
$$\xymatrix{\tV\ar@^{(->}[r]\ar[d]_{\alpha}&\bK\ar[d]^{\beta}\\V_{1}\ar@^{(->}[r]&\bM_{1}}$$
where $\bM_{1}$ is the first blow-up in the Kirwan's process and the vertical maps are blow-ups. A holomorphic map
$$\nu:\tV-\alpha^{-1}(\Delta\cap\tsig/\!/\SL(2))\to\bS$$
was already constructed.

Let $x$ be any point in $\Delta\cap\tsig/\!/\SL(2)$. It follows from (\ref{intersection of D_1^(1) and tsig}) that $x$ is represented by the orbit of $\left[\begin{matrix}a^{0}+p^{0}&0\\0&-a^{0}-p^{0}\end{matrix}\right]$ for some $a^{0}\in H^1(\cO_X)$ and $p^{0}\in H^0(K_X)$. The normal space $Y$ to this orbit is isomorphic to $\cc^{2g}\times\oPsi^{-1}(0)$ where $\cc^{2g}$ is the tangent space of the blow-up $\widetilde{H^1(\cO_X)\oplus H^0(K_X)}=bl_0(H^1(\cO_X)\oplus H^0(K_X))$.

For a neighborhood $Y_1$ of $0$ in $Y$ which is holomorphically embedded into $U_1$, $(\cF_1,\varphi_{\cF_1})|_{X\times Y_{1}}$ has transition matrices of the form
\begin{equation}\label{transition 1}
\left(\begin{matrix}1+z_{1}(a_{ij}^{0}+a_{ij})&z_{1}b_{ij}\\z_{1}c_{ij}&1-z_{1}(a_{ij}^{0}+a_{ij})\end{matrix}\right)
\end{equation}
and local Higgs fields of the form
\begin{equation}\label{Higgs field 1}
\left(\begin{matrix}z_{1}(p_{i}^{0}+p_{i})&z_{1}q_{i}\\z_{1}r_{i}&-z_{1}(p_{i}^{0}+p_{i})\end{matrix}\right)
\end{equation}
where $(a,p)=(\{a_{ij}\},\{p_{i}\})$, $(b,q)=(\{b_{ij}\},\{q_{i}\})$, $(c,r)=(\{c_{ij}\},\{r_{i}\})$ are classes in $H^1(\cO_X)\oplus H^0(K_X)$, not parallel to $(a^0,p^0)$ if they are nonzero and $z_1$ is the coordinate for the normal direction of $\pp(H^1(\cO_X)\oplus H^0(K_X))$ in $\widetilde{H^1(\cO_X)\oplus H^0(K_X)}$.

Since Luna's \'etale slice theorem tells us that a neighborhood of the vertex of the cone $Y/\!/\cc^{\times}$ is analytically equivalent to a neighborhood of $x$ in $V_1$, the image of $\tY$ in $\tV$ is biholomorphic to a neighborhood of $\alpha^{-1}(x)$ where $\tY$ be the strict transform of $Y_1$ in $\tU$.

Following the same argument as in \cite[4.5]{KL}, we can construct a family of rank $4$ Higgs bundles on $X$ parametrized by $\tY$ satisfying the conditions of Proposition \ref{key of construction of morphism}. Then $\nu$ can extends to $\alpha^{-1}(x)$.

\subsection{Contractions} Let $(\cc^{2g},\omega)$ be a symplectic vector space and let $\Gr^{\omega}(k,2g)$ be the Grassmannian of $k$-dimensional subspaces of $\cc^{2g}$ such that the restriction of $\omega$ to the subspace is zero. Let $\cA$ (resp. $\cB$) be the tautological rank $2$ (resp. rank $3$) bundle over the Grassmannian $\Gr^{\omega}(2,2g)$ (resp. $\Gr^{\omega}(3,2g)$). For $B\in\Gr^{\omega}(3,2g)$, let
$$\bC\bC(B):=\text{the closure of }\{(C,D)\in\pp(S^{2}B)\times\pp(S^{2}B^{\vee})|C\text{ and }D\text{ are smooth}$$
$$\text{conics dual to each other}\},$$
which is called the variety of complete conics. Then we have two projections
$$\xymatrix{\pp(S^{2}B)&\bC\bC(B)\ar[l]_{\Phi_{B}}\ar[r]^{\Phi_{B}^{\vee}}&\pp(S^{2}B^{\vee})}$$
where each of $\Phi_{B}$, $\Phi_{\breve{B}}$ is the blow-up of each of $\pp(S^{2}B)$, $\pp(S^{2}B^{\vee})$ along the locus of rank $1$ conics. Let $\bC\bC(\cB)$ be the tautological family of complete conics over $\Gr^{\omega}(3,2g)$. Then the restrictions of the exceptional divisors and their intersections over a point in $\zz_2^{2g}$ can be described as follows. (See also Proposition 4.2 in \cite{KY})

\begin{prop}\label{description of exc div}
\begin{enumerate}
\item[(1)] $\tD^{(1)}\cong\bC\bC(\cB)$, that is, $\tD^{(1)}$ is the blow-up of $\pp(S^{2}\cB)$ along the locus of rank $1$ conics. Hence $\tD^{(1)}$, that is, $\tD^{(1)}$ is a $\hat{\pp}^5$-bundle over $\Gr^\omega(3,2g)$ where $\hat{\pp}^5$ be the blow-up of $\pp^5$ (projectivization of the space of $3\times 3$ symmetric matrices) along $\pp^2$ (the locus of rank $1$ matrices).
\item[(2)] $\tD^{(3)}$ is a $\pp^{2g-4}$-bundle over $\pp(S^{2}\cA)$. Hence $\tD^{(3)}$ is a $\pp^{2g-4}$-bundle over a $\pp^2$-bundle over $\Gr^{\omega}(2,2g)$.
\item[(3)] $\tD^{(1)}\cap\tD^{(3)}$ is the proper transform of $\pp(S^{2}\cB)_2$(the locus of rank $\leq 2$ matrices) in the blow-up $\bC\bC(\cB)\to\pp(S^{2}\cB)$, that is, the exceptional divisor of the blow-up $\bC\bC(\cB)\to\pp(S^{2}\cB^{\vee})$. Hence $\tD^{(1)}\cap\tD^{(3)}$ is a $\pp^2$-bundle over a $\pp^2$-bundle over $\Gr^\omega(3,2g)$.
\item[(4)] $\tD^{(1)}\cap\tD^{(2)}\cap\tD^{(3)}$ is a $\pp^1$-bundle over a $\pp^2$-bundle over $\Gr^\omega(3,2g)$.
\item[(5)] $\tD^{(1)}\cap\tD^{(2)}$ is the exceptional divisor of the blow-up $\bC\bC(\cB)\to\pp(S^{2}\cB)$. Hence $\tD^{(1)}\cap\tD^{(2)}$ is a $\pp^2$-bundle over a $\pp^2$-bundle over $\Gr^\omega(3,2g)$.
\end{enumerate}
\end{prop}
We next investigate some rational curves that will be contracted. Define
$$\s:=\text{ class of a line in a }\pp^2\text{-fiber of }\Phi_{B}^{\vee}$$
$$\e:=\text{ class of a line in a }\pp^2\text{-fiber of }\Phi_{B}$$
$$\g:=\text{ class of }\{\Phi_{B_t}^{-1}([q_t])\}_{t\in\Lambda}$$
in $N_{1}(\tD^{(1)})$ (the group of numerical equivalence classes of $1$-cycles) where $[A]\in\Gr^\omega(2,2g)$, $q\in S^{2}A$, $A^{\perp}$ is the orthogonal complement of $A$ with respect to $\omega$, a line $\Lambda\subset\pp(A^{\perp}/A)$ are fixed, $B_t\subset\cc^{2g}$ is the $3$-dimensional subspaces containing $A$ for $t\in\Lambda$ (that is, $[B_t]\in\Gr^\omega(3,2g)$) and $q_t\in S^{2}B_t$ is the image of $q$ under the inclusion $S^{2}A\hookrightarrow S^{2}B_t$.

Let $\hs=\iota_*\s$, $\he=\iota_*\e$ and $\hg=\iota_*\g$ where $\iota$ is the inclusion $\tD^{(1)}\hookrightarrow\bK$. The following proposition tells us that $\bK$ can be blown-down twice.
\begin{prop}\label{contraction}
\begin{enumerate}
\item[(1)] $\rr^{+}\hs$ is a $\omega_{\bK}$-negative extremal ray. The contraction $\bK_{\s}$ of the ray $\rr^{+}\hs$ is a smooth quasi-projective desingularization of $\bM_2$. In fact, this is the contraction of the fiber direction of $\pp(S^{2}\cA)\to\Gr^\omega(2,2g)$ in $\tD^{(3)}$ and is also a blow-down map.
\item[(2)] For the image $\bar{\e}$ of $\he$ in $N_1(\bK_{\s})$, $\rr^{+}\bar{\e}$ is $\omega_{\bK_{\s}}$-negative extremal ray and its contraction $\bK_{\e}$ is a smooth quasi-projective desingularization of $\bM_2$. This is the contraction of the fiber direction of $\pp(S^{2}\cB^{\vee})\to\Gr^\omega(3,2g)$ and is also a blow-down map.
\end{enumerate}
\end{prop}
\begin{proof}
\begin{enumerate}
The proofs are identical to those of (3.0.2)-(3.0.4) in \cite{OGr97}.
\end{enumerate}
\end{proof}

\subsection{A stratification of $\bS$}
In this subsection, we provide a geometric description for the canonical subschemes of $\bS$, which is useful for the proof of Theorem \ref{main}. The arguments are almost same as those of \cite{Bal} except for the necessary modifications. Since $\bS$ is a fine moduli scheme, we have the universal family $(\cE_{univ},\phi_{univ},s_{univ})$ of rank $4$ and degree $0$ on $X\times\bS$. Consider the canonical family of specializations $B=p_{\bS *}\EEnd((\cE_{univ},\phi_{univ}))$ of $M(2)$ parametrized by $\bS$ and the universal family of specializations $W$ of $M(2)$ parametrized by $\cA_2$, where $p_{\bS}:X\times\bS\to\bS$ is the projection onto the second component. Then $B$ and $W$ give generalized conic bundles $P$ on $\bS$ and $Q$ on $\cA_2$ respectively. (See \cite[Remark 3]{Bal} for more details)

\begin{prop}[\cite{Bal}, Proposition 2]\label{base change of conic bundles}
$P=\pi_{\text{sp}}^*Q$ for the morphism $\pi_{\text{sp}}:\bS\to\cA_2$ given by $(E,\phi,s)\mapsto\End((E,\phi))$.
\end{prop}
\begin{proof}

This is an immediate consequence of the definitions of $\pi_{sp}$, $B$ and $W$.
\end{proof}

To define the canonical subschemes of $\cA_2$, we need the concept of \textit{generalized conic bundle} in the sense of \cite{Bal} as follows.

\begin{Def}[\cite{Bal}, Definition 1]\label{gen conic bdl}
Let $Y$ be a variety. A \emph{generalized conic bundle} $\cC$ on $Y$ is defined by
\begin{enumerate}
\item[(a)] a vector bundle $V$ on $Y$ of rank $3$ and

\item[(b)] a closed subscheme $\cC$ of $\pp(V)$ over $Y$, such that, given $y\in Y$, there exists a neighborhood $U$ of $y$, where $\cC\cap p^{-1}(U)$ is defined by $q=0$, $q\in H^0(p^{-1}(U),H^2)$, $H$ being the tautological line bundle for $p:\pp(V)\to Y$, that is, $p_{*}H\cong V^{\vee}$.
\end{enumerate}
\end{Def}

Since $\cC$ is an effective Cartier divisor by Definition \ref{gen conic bdl}, it corresponds to a section of a line bundle $L_{\cC}$ on $\pp(V)$. Since $L_{\cC}$ and $H^2$ coincide locally over $Y$, it follows from the ``see-saw" theorem (cf. Mumford's \textit{Abelian varieties}) that there exists a line bundle $L$ on $Y$ such that $L_{\cC}=H^{2}\otimes p^*(L)$. Since $p_{*}L_{\cC}=p_*(H^2)\otimes L=S^2(V^*)\otimes L$, $q$ of the condition (b) in Definition \ref{gen conic bdl} can be replaced by an element $q_{L}\in H^0(S^2(V^*)\otimes L)$, that is, a \textit{quadratic form} $q_{L}:V\to L$. The \textit{discriminant} $\Delta$ of $q_{L}$ is defined as a section of $L^3\otimes(\wedge^3(V^{\vee}))^2$. The \textit{degeneracy locus} of $\cC$ is given by $\Delta=0$.

For $i=1,2,3,$ set
$$Y_{i}=\{y\in Y|q_{L}|_{V|_{y}}\text{ has rank }\leq3-i\}.$$
Then we have $Y\supset Y_{1}\supset Y_{2}\supset Y_{3}$. We call $Y_{i}$ the \textit{canonical subschemes associated to the degenerate loci of the conic bundle} $\cC$ on $Y$.

We have the canonical subschemes
$$\bS\supset \bS_{1}\supset \bS_{2}\supset \bS_{3}$$
and
$$\cA_2\supset(\cA_2)_{1}\supset (\cA_2)_{2}\supset (\cA_2)_{3}$$
associated to the degeneracy locus of $P$ and $Q$ respectively. Then, by Proposition \ref{base change of conic bundles}, $\pi_{sp}$ maps $\bS-\bS_{2}$ into $\cA_{2}-(\cA_{2})_{2}$ such that $\bS_{1}-\bS_{2}\to(\cA_{2})_{1}-(\cA_{2})_{2}$ and $\bS-\bS_{1}\to\cA_{2}-(\cA_{2})_{1}$.

By Theorem 1 in \cite{Se77}, we have $\cA_{2}\simeq\Phi\times\Lambda$, where $\Lambda$ is the $3$-dimensional affine space and $\Phi$ is the $6$-dimensional affine space whose points are identified with the set of quadratic forms on a fixed $3$-dimensional vector space or the set of algebras of the form $C_{q}^{+}$, the even degree elements of the Clifford algebra associated to the quadratic form $q$. Hence we have
$$(\cA_{2})_{i}=\{q\in\Phi|\text{rank }q\leq 3-i\}\times\cc^3\text{ for }i=1,2,3,$$
and
$$\cA_{2}-(\cA_{2})_{1}=\{q\in\Phi|C_{q}^{+}\simeq M(2)\}\times\cc^3=\{y\in\cA_{2}|W_{y}\simeq M(2)\}.$$

\begin{prop}\label{geometry of canonical subschemes}
The subsets $\bS-\bS_{2}$ and $\bS_{1}-\bS_{2}$ of $\bS$ are precisely $\bS-\pi_{\bS}^{-1}(\zz_2^{2g})$ and $\pi_{\bS}^{-1}(T^*J/\zz_{2}-\zz_2^{2g})$ respectively. In particular,  $\bS-\bS_{1}=\pi_{\bS}^{-1}(\bM_2^s)$.
\end{prop}
\begin{proof}
The proof is similar to that of Proposition 3 in \cite{Bal}. By the arguments before this proposition, it suffices to show that the subsets $\pi_{\bS}^{-1}(\bM_2^s)$ and $\pi_{\bS}^{-1}(T^*J/\zz_{2}-\zz_2^{2g})$ of $\bS$ are mapped by $\pi_{sp}$ into the subsets $\cA_2-(\cA_2)_1$ and $(\cA_2)_{1}-(\cA_2)_{2}$ of $\cA_2$ respectively. We already know that $(E,\phi,s)\in\pi_{\bS}^{-1}(\bM_2^s)$ if and only if $\End((E,\phi))=M(2)$, which proves that $\pi_{\bS}^{-1}(\bM_2^s)$ maps to $\cA_{2}-(\cA_{2})_1$.

Next, we claim that for $(E,\phi,s)\in\pi_{\bS}^{-1}(T^*J/\zz_{2}-\zz_2^{2g})$, $\End((E,\phi))$ has the same defining relations as that of $C_{q}^{+}$ for a quadratic form $q$ of rank $2$ on a $3$-dimensional  vector space. It follows from the definition of $\pi_\bS$ that the endomorphism algebras of any two points in a fibre $\pi_{\bS}^{-1}((L,\psi)\oplus(L^{-1},-\psi))$ are isomorphic. So we can choose a point $(E,\phi,s)\in\pi_{\bS}^{-1}((L,\psi)\oplus(L^{-1},-\psi))$, where
$$(E,\phi)=(V,\phi_{V})\oplus(W,\phi_{W}),$$
$$(V,\phi_{V})\in H^1(L^2)\oplus H^0(L^2 K_X),\quad(W,\phi_{W})\in H^1(L^{-2})\oplus H^0(L^{-2}K_X)$$
and
$$(L,\psi)\in T^*J-\zz_2^{2g},$$
that is,
\begin{equation}\label{extension1}0\to(L,\psi)\to(V,\phi_{V})\to(L^{-1},-\psi)\to0\end{equation}
and
\begin{equation}\label{extension2}0\to(L^{-1},-\psi)\to(W,\phi_{W})\to(L,\psi)\to0\end{equation}
It is obvious that this kinds of points belongs to $\pi_{\bS}^{-1}(T^*J/\zz_{2}-\zz_2^{2g})$. By (\ref{extension1}) and (\ref{extension2}), we can see that $\End((V,\phi_{V})\oplus(W,\phi_{W}))$ has four generators, which can be written as the following block matrices
$$x=\left(\begin{matrix}0&0\\\gamma_{2}&0\end{matrix}\right),
w=\left(\begin{matrix}0&\gamma_{1}\\0&0\end{matrix}\right), u=\left(\begin{matrix}I&0\\0&0\end{matrix}\right)\text{ and } v=\left(\begin{matrix}0&0\\0&I\end{matrix}\right)$$
where $I$ is the $2\times 2$ identity matrix and $\gamma_{1},\gamma_{2}$ are the $2\times 2$ matrices that come from identification of the rank $1$ Higgs bundles in (\ref{extension1}) and (\ref{extension2}).

Here we need to show that $\End((V,\phi_V))=\cc\cdot\id_V$ precisely. First of all we show that $\End(V)=\cc\cdot\id_V$. Consider an arbitrary nontrivial morphism $\phi:V\to V$ in $\End(V)$. Let
$$\xymatrix{0\ar[r]&L\ar[r]^{i}&V\ar[r]^{j}&L^{-1}\ar[r]&0}$$
be the nonsplit extension to define $V$. Then the composition $$\xymatrix{L\ar[r]^{i}&V\ar[r]^{\phi}&V\ar[r]^{j}&L^{-1}}$$
should be zero. If not, then this composition is an isomorphism, which is a contradiction. Then $\phi(L)\subset L$. Let $\psi:=\phi-\lambda\cdot\id_V$. Since $\phi|_{L}=\lambda\cdot\id_L$ for some $\lambda\in\cc^\times$, $\psi|_{L}=0$, that is, $\psi\circ i=0$, which implies that $\psi:V\to V$ factors through $\xymatrix{V\ar[r]^{j}&L^{-1}}$ and there exists a unique morphism $\psi':L^{-1}\to V$ such that $\psi'\circ j=\psi$. If the composition $\xymatrix{L^{-1}\ar[r]^{\psi'}&V\ar[r]^{j}&L^{-1}}$ is nonzero, then $j\circ\psi'=\eta\cdot\id_{L^{-1}}$ for some $\eta\in\cc^\times$. Then $\frac{1}{\eta}\psi':L^{-1}\to V$ is a splitting of the following nonplit extension
$$\xymatrix{0\ar[r]&L\ar[r]^{i}&V\ar[r]^{j}&L^{-1}\ar[r]&0},$$
which is a contradiction. So the comosition $\xymatrix{L^{-1}\ar[r]^{\psi'}&V\ar[r]^{j}&L^{-1}}$ is zero and then there exists a unique morphism $h:L^{-1}\to L$ such that $i\circ h=\psi'$. Since $L\ncong L^{-1}$, $h=0$, that is, $\psi'=0$. Hence $\psi=0$, that is, $\psi=\lambda\cdot\id_V$.\\
By Theorem \ref{def sp of Higgs}, $\End((V,\psi_V))=\bH^0\bfHom((V,\phi_V),(V,\phi_V))$. We know that
$$\bH^0\bfHom((V,\phi_V),(V,\phi_V))=\ker(\xymatrix{\End(V)\ar[r]^{[\quad,\phi_V]\quad}&\Hom(V,V\otimes K_X)}).$$
Since $\End(V)=\cc\cdot\id_V$, we have $[\quad,\phi_V]=0$, that is, $\End((V,\psi_V))=\End(V)=\cc\cdot\id_V$.\\
Similarly, we can show that $\End((W,\phi_W))=\cc\cdot\id_W$.

Coming back to the main proof, the defining relations of $\End((V,\phi_{V})\oplus(W,\phi_{W}))$ can be written as
$$u^2=u, v^2=v, uv=0, u+v=I,$$
\begin{equation}\label{def rel1}w^2=x^2=wx=0, uw=w, wu=0,\end{equation}
$$ux=0, xu=x, vw=0, wv=w,$$
$$vx=x, xv=0.$$
For a quadratic form $q$ of rank $2$ on a $3$-dimensional vector space over $\cc$, we can see that $C_q^+$ is a $4$-dimensional $\cc$-algebra with
$$C_{q}^{+}=\cc+\cc\alpha+\cc\beta+\cc\gamma\text{ such that}$$
$$\alpha^{2}=-1, \alpha\beta=-\gamma, \alpha\gamma=\beta, \beta\alpha=\gamma, \gamma\alpha=-\beta.$$
When we set $a=\frac{1}{2}(1+i\alpha), b=\frac{1}{2}(1-i\alpha), c=i\beta+\gamma\text{ and }d=i\beta-\gamma$  for $i=\sqrt{-1}\in\cc$, $a,b,c,d$ are new generators of $C_q^+$ with the following defining relations
$$a^2=a, b^2=b, ab=0, a+b=1,$$
\begin{equation}\label{def rel2}c^2=d^2=cd=0, ac=c, ca=0,\end{equation}
$$ad=0, da=d, bc=0, cb=c,$$
$$bd=d, db=0.$$
Both (\ref{def rel1}) and (\ref{def rel2}) are the same.
\end{proof}

\begin{Cor}\label{geometry1 of 1-2 canonical subscheme}
$\bS_{1}-\bS_{2}\cong\tD^{(2)}-(\tD^{(1)}\cup\tD^{(3)})$, that is, it is a free $\zz_{2}$-quotient of a $I_{2g-3}$-bundle on $T^*J-\zz_2^{2g}$, where $I_{2g-3}$ is the incidence variety given by $I_{2g-3}=\{(p,H)\in\pp^{2g-3}\times(\pp^{2g-3})^{\ast}\;|\;p\in H\}$.
\end{Cor}
\begin{proof}
The proof is similar to that of \cite[Corollary 1]{Bal} except for using (\ref{def rel1}) and Propositon \ref{surjective}.

We claim that if $(E,\phi,s)\in\bS_{1}-\bS_{2}=\pi_{\bS}^{-1}(T^*J/\zz_{2}-\zz_2^{2g})$ and $\pi_{\bS}((E,\phi,s))=(L,\psi)\oplus(L^{-1},-\psi)\in T^*J/\zz_{2}-\zz_2^{2g}$ then $(E,\phi)=(V,\phi_{V})\oplus(W,\phi_{W})$ for some extension classes
$$[(V,\phi_{V})]\in\pp\bH^{1}\bfHom((L,\psi),(L^{-1},-\psi))=\pp(H^{1}(L^{-2})\oplus H^{0}(L^{-2}K_{X}))=\pp^{2g-3}$$
and
$$[(W,\phi_{W})]\in\pp\bH^{1}\bfHom((L^{-1},-\psi),(L,\psi))=\pp(H^{1}(L^{2})\oplus H^{0}(L^{2}K_{X}))=\pp^{2g-3}.$$

Let $(E,\phi,s)\in\pi_{\bS}^{-1}(T^*J/\zz_{2}-\zz_2^{2g})$ such that $\pi_{\bS}((E,\phi,s))=(L,\psi)\oplus(L^{-1},-\psi)\in T^*J/\zz_{2}-\zz_2^{2g}$. Then $\End((E,\phi))$ has four generators $x,w,u,v$ with defining relations (\ref{def rel1}) as in the proof of Proposition \ref{geometry of canonical subschemes}. Let $(V,\phi_{V})=\ker u$. Then $(V,\phi_{V})$ is a Higgs subbundle of $(E,\phi)$ and we have an exact sequence
$$0\to(V,\phi_{V})\to(E,\phi)\to(W,\phi_{W})\to0$$
such that $(W,\phi_{W})=\ker v$. Thus this exact sequence splits and then $(E,\phi)=(V,\phi_{V})\oplus(W,\phi_{W})$.

By Propositon \ref{surjective}, $(V,\phi_{V})$ and $(W,\phi_{W})$ cannot be $(L,\psi)\oplus(L,\psi)$ or $(L^{-1},-\psi)\oplus(L^{-1},-\psi)$ and we also rule out $(V,\phi_{V})=(L,\psi)\oplus(L^{-1},-\psi)$ and $(W,\phi_{W})=(L^{-1},-\psi)\oplus(L,\psi)$. Hence
$$[(V,\phi_{V})]\in\pp\bH^{1}\bfHom((L,\psi),(L^{-1},-\psi))\text{ and }[(W,\phi_{W})]\in\pp\bH^{1}\bfHom((L^{-1},-\psi),(L,\psi)).$$

We next claim that if we choose $(c,d,e,f)\in[H^{0}(L^{-2}K_{X})\oplus H^{0}(L^{2}K_{X})]\oplus[H^{1}(L^{2})\oplus H^{1}(L^{-2})]$ such that $[(V,\phi_{V})]=[(f,c)]$ and $[(W,\phi_{W})]=[(e,d)]$ then $(c,d,e,f)\in\oPsi^{-1}(0)$.

Assume that $(c,d,e,f)\not\in\oPsi^{-1}(0)$. Let $(a,b)\in H^0(K_{X})\oplus H^{1}(\cO_{X})$ be classes obtained from transition of $L$ and local Higgs field of $\psi$. Following the arguments of the proof of Lemma \ref{useful lemma for the assumption}, $(V,\phi_{V})$ and $(W,\phi_{W})$ are obtained by an elementary modification starting with a family $(\cF,\varphi_{\cF})$ of semistable Higgs bundles over $X$ parametrized by
$$C=\{(a,b,zc,zd,ze,zf)\in H^0(K_{X})\oplus H^{1}(\cO_{X})\oplus[H^{0}(L^{-2}K_{X})\oplus H^{0}(L^{2}K_{X})]\oplus[H^{1}(L^{2})\oplus H^{1}(L^{-2})]\;|\;z\in\cc\}.$$
Precisely, $(V,\phi_{V})=(\cF',\varphi_{\cF'})|_{z=0}$ and $(W,\phi_{W})=(\cF'',\varphi_{\cF''})|_{z=0}$, where $(\cF',\varphi_{\cF'})$ is the kernel of $(\cF,\varphi_{\cF})\to(\cF,\varphi_{\cF})|_{z=0}\cong(L^{-1},-\psi)\oplus(L^{-1},-\psi)\to(L,\psi)$ and $(\cF'',\varphi_{\cF''})$ is the kernel of $(\cF,\varphi_{\cF})\to(\cF,\varphi_{\cF})|_{z=0}\cong(L^{-1},-\psi)\oplus(L^{-1},-\psi)\to(L^{-1},-\psi)$.

But since $(c,d,e,f)\not\in\oPsi^{-1}(0)$, $(zc,zd,ze,zf)\not\in\oPsi^{-1}(0)$ for all $z\in\cc$ and then $(\cF,\varphi_{\cF})|_{X\times(a,b,zc,zd,ze,zf)}$ is not a Higgs bundle for all $z\in\cc$. Then $(V,\phi_{V})$ and $(W,\phi_{W})$ are not Higgs bundles. Indeed, if $(V,\phi_{V})$ and $(W,\phi_{W})$ are Higgs bundles and if we consider the following elementary modifications
$$0\to(\cG',\varphi_{\cG'})\to(\cF',\varphi_{\cF'})\to(L^{-1},-\psi)\to0$$
and
$$0\to(\cG'',\varphi_{\cG''})\to(\cF'',\varphi_{\cF''})\to(L,\psi)\to0,$$
then $(\cG',\varphi_{\cG'})\cong(\cG'',\varphi_{\cG''})\cong(\cF\otimes\cO_{X\times C}(-D),\varphi_{\cF\otimes\cO_{X\times C}(-D)})$, where $D=\{(x,(a,b,zc,zd,ze,zf))\in X\times C\;|\;z=0\}$. Since $(\cF\otimes\cO_{X\times C}(-D),\varphi_{\cF\otimes\cO_{X\times C}(-D)})|_{z=0}$ is a Higgs bundle and $\cO_{X\times C}(-D)$ is a line bundle, $(V,\phi_{V})$ and $(W,\phi_{W})$ are Higgs bundles. Thus we get a contradiction. Hence $(c,d,e,f)\in\oPsi^{-1}(0)$.

\end{proof}

\subsection{Factorization of $\rho$}

In this subsection, we show that $\rho$ factors through $\bK_\e$. Let $\bD_1$, $\bD_2$ and $\bD_3$ be three exceptional divisors on $\bK$ coming from the three blow-ups. Consider the point $0=[(\cO_X,0)\oplus(\cO_X,0)]\in\zz_2^{2g}$. Let $\tD^{(1)}$, $\tD^{(2)}$ and $\tD^{(3)}$ be the restrictions of $\bD_1$, $\bD_2$ and $\bD_3$ to a neighborhood of the preimage of $0$ in $\bK$ respectively. We need the following two lemmas for the proof.

\begin{Lem}\label{indep1}
\begin{enumerate}
\item The isomorphism classes of Higgs bundles given by (\ref{transition for cE_1|_{z=0}}) and (\ref{Higgs field of cE_1|_{z=0}}) are independent of the choice of $(a,d)$, $(b,e)$.
\item The isomorphism classes of Higgs bundles given by (\ref{transition for tcE 1}) and (\ref{Higgs field of tcE 1}) are independent of the choice of $(g,h)$, $(b,e)$, $(c,f)$.
\end{enumerate}
\end{Lem}
\begin{proof} The proofs of (1) and (2) are identical. We prove only (1). Let $\{V_i\}_{i}$ be an open cover of $X$ and let $(E,\phi)$ and $(E',\phi')$ be Higgs bundles given by (\ref{transition for cE_1|_{z=0}}) on $V_{i}\cap V_{j}$ and (\ref{Higgs field of cE_1|_{z=0}}) on $V_i$. Assume that an isomorphism between $(E,\phi)$ and $(E',\phi')$ is given by invertible matrices
$$\left(
\begin{matrix}\eta_{i,11}&\eta_{i,12}&\eta_{i,13}&\eta_{i,14}\\
\eta_{i,21}&\eta_{i,22}&\eta_{i,23}&\eta_{i,24}\\
\eta_{i,31}&\eta_{i,32}&\eta_{i,33}&\eta_{i,34}\\
\eta_{i,41}&\eta_{i,42}&\eta_{i,43}&\eta_{i,44}\end{matrix}\right)$$
on $V_i$.
Then
$$
\left(
\begin{matrix}\eta_{j,11}&\eta_{j,12}&\eta_{j,13}&\eta_{j,14}\\
\eta_{j,21}&\eta_{j,22}&\eta_{j,23}&\eta_{j,24}\\
\eta_{j,31}&\eta_{j,32}&\eta_{j,33}&\eta_{j,34}\\
\eta_{j,41}&\eta_{j,42}&\eta_{j,43}&\eta_{j,44}\end{matrix}\right)\left(
\begin{matrix}1&0&0&0\\
0&1&0&0\\
0&0&1&0\\
-2\alpha_{ij}|_{z=0}&-\beta_{ij}|_{z=0}&\gamma_{ij}|_{z=0}&1\end{matrix}\right)$$
$$=\left(
\begin{matrix}1&0&0&0\\
0&1&0&0\\
0&0&1&0\\
-2\alpha'_{ij}|_{z=0}&-\beta'_{ij}|_{z=0}&\gamma'_{ij}|_{z=0}&1\end{matrix}\right)\left(
\begin{matrix}\eta_{i,11}&\eta_{i,12}&\eta_{i,13}&\eta_{i,14}\\
\eta_{i,21}&\eta_{i,22}&\eta_{i,23}&\eta_{i,24}\\
\eta_{i,31}&\eta_{i,32}&\eta_{i,33}&\eta_{i,34}\\
\eta_{i,41}&\eta_{i,42}&\eta_{i,43}&\eta_{i,44}\end{matrix}\right)
$$
and
$$
\left(
\begin{matrix}\eta_{i,11}&\eta_{i,12}&\eta_{i,13}&\eta_{i,14}\\
\eta_{i,21}&\eta_{i,22}&\eta_{i,23}&\eta_{i,24}\\
\eta_{i,31}&\eta_{i,32}&\eta_{i,33}&\eta_{i,34}\\
\eta_{i,41}&\eta_{i,42}&\eta_{i,43}&\eta_{i,44}\end{matrix}\right)\left(
\begin{matrix}0&0&0&0\\
0&0&0&0\\
0&0&0&0\\
-2p_i|_{z=0}&-q_i|_{z=0}&r_i|_{z=0}&0\end{matrix}\right)
$$
$$=\left(
\begin{matrix}0&0&0&0\\
0&0&0&0\\
0&0&0&0\\
-2p'_i|_{z=0}&-q'_i|_{z=0}&r'_i|_{z=0}&0\end{matrix}\right)\left(
\begin{matrix}\eta_{i,11}&\eta_{i,12}&\eta_{i,13}&\eta_{i,14}\\
\eta_{i,21}&\eta_{i,22}&\eta_{i,23}&\eta_{i,24}\\
\eta_{i,31}&\eta_{i,32}&\eta_{i,33}&\eta_{i,34}\\
\eta_{i,41}&\eta_{i,42}&\eta_{i,43}&\eta_{i,44}\end{matrix}\right).$$
Equating each rows, we have
$$\eta_{i,11}=\eta_{j,11}=\text{constant},\eta_{i,12}=\eta_{j,12}=\text{constant},\eta_{i,13}=\eta_{j,13}=\text{constant},\eta_{i,14}=\eta_{j,14}=0,$$
$$\eta_{i,21}=\eta_{j,21}=\text{constant},\eta_{i,22}=\eta_{j,22}=\text{constant},\eta_{i,23}=\eta_{j,23}=\text{constant},\eta_{i,24}=\eta_{j,24}=0,$$
$$\eta_{i,31}=\eta_{j,31}=\text{constant},\eta_{i,32}=\eta_{j,32}=\text{constant},\eta_{i,33}=\eta_{j,33}=\text{constant},\eta_{i,34}=\eta_{j,34}=0,$$
$$\eta_{i,44}=\eta_{j,44}=\text{constant}\ne0,$$
and
$$
\begin{matrix}
-2\eta_{i,44}(a,d)=-2\eta_{i,11}(a',d')-\eta_{i,21}(b',e')+\eta_{i,31}(c',f')\\
-\eta_{i,44}(b,e)=-2\eta_{i,12}(a',d')-\eta_{i,22}(b',e')+\eta_{i,32}(c',f')\\
\eta_{i,44}(c,f)=-2\eta_{i,13}(a',d')-\eta_{i,23}(b',e')+\eta_{i,33}(c',f'),
\end{matrix}
$$
where
$$(\{\alpha_{ij}|_{z=0}\},\{p_i|_{z=0}\}), (\{\alpha'_{ij}|_{z=0}\},\{p'_i|_{z=0}\}), (\{\beta_{ij}|_{z=0}\},\{q_i|_{z=0}\}),$$
$$(\{\beta'_{ij}|_{z=0}\},\{q'_i|_{z=0}\}), (\{\gamma_{ij}|_{z=0}\},\{r_i|_{z=0}\})\text{ and }(\{\gamma'_{ij}|_{z=0}\},\{r'_i|_{z=0}\})$$
represents classes $(a,d)$, $(a',d')$, $(b,e)$, $(b',e')$, $(c,f)$ and $(c',f')$ in $H^1(\cO_X)\oplus H^0(K_X)$ respectively.

Since
$$\left(\begin{matrix}-2\eta_{i,11}&-\eta_{i,21}&\eta_{i,31}\\
-2\eta_{i,12}&-\eta_{i,22}&\eta_{i,32}\\
-2\eta_{i,13}&-\eta_{i,23}&\eta_{i,33}\end{matrix}\right)
=\left(\begin{matrix}\eta_{i,11}&\eta_{i,21}&\eta_{i,31}\\
\eta_{i,12}&\eta_{i,22}&\eta_{i,32}\\
\eta_{i,13}&\eta_{i,23}&\eta_{i,33}\end{matrix}\right)
\left(\begin{matrix}-2&0&0\\
0&-1&0\\
0&0&1\end{matrix}\right)$$
is invertible, $\Sp\{(a,d),(b,e),(c,f)\}=\Sp\{(a',d'),(b',e'),(c',f')\}$.
\end{proof}

\begin{Lem}\label{indep2}
The isomorphism classes of Higgs bundles given by (\ref{transition for tcE 2}) and (\ref{Higgs field of tcE 2}) are independent of the choice of $(b,e)$, $(c,f)$.
\end{Lem}
\begin{proof} Let $\{V_i\}_{i}$ be an open cover of $X$ and let $(E,\phi)$ and $(E',\phi')$ be Higgs bundles given by (\ref{transition for tcE 2}) on $V_{i}\cap V_{j}$ and (\ref{Higgs field of tcE 2}) on $V_i$. Assume that an isomorphism between $(E,\phi)$ and $(E',\phi')$ is given by invertible matrices
$$\left(
\begin{matrix}\eta_{i,11}&\eta_{i,12}&\eta_{i,13}&\eta_{i,14}\\
\eta_{i,21}&\eta_{i,22}&\eta_{i,23}&\eta_{i,24}\\
\eta_{i,31}&\eta_{i,32}&\eta_{i,33}&\eta_{i,34}\\
\eta_{i,41}&\eta_{i,42}&\eta_{i,43}&\eta_{i,44}\end{matrix}\right)$$
on $V_i$.
Then
$$
\left(
\begin{matrix}\eta_{j,11}&\eta_{j,12}&\eta_{j,13}&\eta_{j,14}\\
\eta_{j,21}&\eta_{j,22}&\eta_{j,23}&\eta_{j,24}\\
\eta_{j,31}&\eta_{j,32}&\eta_{j,33}&\eta_{j,34}\\
\eta_{j,41}&\eta_{j,42}&\eta_{j,43}&\eta_{j,44}\end{matrix}\right)\left(
\begin{matrix}1&0&0&0\\
\gamma_{ij}|_{z=0}&1&0&0\\
\beta_{ij}|_{z=0}&0&1&0\\
-2\delta_{ij}|_{z=0}&-\beta_{ij}|_{z=0}&\gamma_{ij}|_{z=0}&1\end{matrix}\right)$$
$$=\left(
\begin{matrix}1&0&0&0\\
\gamma'_{ij}|_{z=0}&1&0&0\\
\beta'_{ij}|_{z=0}&0&1&0\\
-2\delta'_{ij}|_{z=0}&-\beta'_{ij}|_{z=0}&\gamma'_{ij}|_{z=0}&1\end{matrix}\right)\left(
\begin{matrix}\eta_{i,11}&\eta_{i,12}&\eta_{i,13}&\eta_{i,14}\\
\eta_{i,21}&\eta_{i,22}&\eta_{i,23}&\eta_{i,24}\\
\eta_{i,31}&\eta_{i,32}&\eta_{i,33}&\eta_{i,34}\\
\eta_{i,41}&\eta_{i,42}&\eta_{i,43}&\eta_{i,44}\end{matrix}\right)
$$
and
$$
\left(
\begin{matrix}\eta_{i,11}&\eta_{i,12}&\eta_{i,13}&\eta_{i,14}\\
\eta_{i,21}&\eta_{i,22}&\eta_{i,23}&\eta_{i,24}\\
\eta_{i,31}&\eta_{i,32}&\eta_{i,33}&\eta_{i,34}\\
\eta_{i,41}&\eta_{i,42}&\eta_{i,43}&\eta_{i,44}\end{matrix}\right)\left(
\begin{matrix}0&0&0&0\\
r_i|_{z=0}&0&0&0\\
q_i|_{z=0}&0&0&0\\
-2s_i|_{z=0}&-q_i|_{z=0}&r_i|_{z=0}&0\end{matrix}\right)
$$
$$=\left(
\begin{matrix}0&0&0&0\\
r'_i|_{z=0}&0&0&0\\
q'_i|_{z=0}&0&0&0\\
-2s'_i|_{z=0}&-q'_i|_{z=0}&r'_i|_{z=0}&0\end{matrix}\right)\left(
\begin{matrix}\eta_{i,11}&\eta_{i,12}&\eta_{i,13}&\eta_{i,14}\\
\eta_{i,21}&\eta_{i,22}&\eta_{i,23}&\eta_{i,24}\\
\eta_{i,31}&\eta_{i,32}&\eta_{i,33}&\eta_{i,34}\\
\eta_{i,41}&\eta_{i,42}&\eta_{i,43}&\eta_{i,44}\end{matrix}\right).$$
Equating each rows, we have
$$\eta_{i,11}=\eta_{j,11}=\text{constant},\eta_{i,12}=\eta_{j,12}=0,\eta_{i,13}=\eta_{j,13}=0,\eta_{i,14}=\eta_{j,14}=0,$$
$$\eta_{i,22}=\eta_{j,22}=\text{constant},\eta_{i,23}=\eta_{j,23}=\text{constant},\eta_{i,24}=\eta_{j,24}=0,$$
$$\eta_{i,32}=\eta_{j,32}=\text{constant},\eta_{i,33}=\eta_{j,33}=\text{constant},\eta_{i,34}=\eta_{j,34}=0,$$
$$\eta_{i,44}=\eta_{j,44}=\text{constant},$$
and
\begin{equation}\label{linear relation1}
\begin{matrix}
-\eta_{i,22}(c,f)-\eta_{i,23}(b,e)+\eta_{i,11}(c',f')=0\\
-\eta_{i,32}(c,f)-\eta_{i,33}(b,e)+\eta_{i,11}(b',e')=0\\
-\eta_{i,22}(b',e')+\eta_{i,32}(c',f')+\eta_{i,44}(b,e)=0\\
-\eta_{i,23}(b',e')+\eta_{i,33}(c',f')-\eta_{i,44}(c,f)=0
\end{matrix}
\end{equation}
where
$$(\{\delta_{ij}|_{z=0}\},\{s_i|_{z=0}\}), (\{\delta'_{ij}|_{z=0}\},\{s'_i|_{z=0}\}), (\{\beta_{ij}|_{z=0}\},\{q_i|_{z=0}\}),$$
$$(\{\beta'_{ij}|_{z=0}\},\{q'_i|_{z=0}\}), (\{\gamma_{ij}|_{z=0}\},\{r_i|_{z=0}\})\text{ and }(\{\gamma'_{ij}|_{z=0}\},\{r'_i|_{z=0}\})$$
represents classes $(g,h)$, $(g',h')$, $(b,e)$, $(b',e')$, $(c,f)$ and $(c',f')$ in $H^1(\cO_X)\oplus H^0(K_X)$ respectively. Hence
$$\Sp\{(b,e),(c,f)\}=\Sp\{(b',e'),(c',f')\}.$$
\end{proof}

\begin{prop}\label{factors through}
$\rho$ factors through $\bK_{\e}$.
\end{prop}
\begin{proof}
Consider the first contraction $f_{\s}:\bK\to\bK_{\s}$. We
show that there is a continuous map $\rho_{\s}:\bK_{\s}\to\bS$
such that $\rho_{\s}\circ f_{\s}=\rho$. It is enough to show that $\rho$ is constant on the fibers
of $f_{\s}$ by Riemann's extension theorem. By Proposition \ref{description of exc div} and Proposition \ref{contraction}, $f_{\s}$ is indeed the contraction of the fibers $\pp^2$ of $\pp(S^{2}\cA)\to\Gr^\omega(2,2g)$ in $\tD^{(3)}$ which forgets the choice of $(b,e)$, $(c,f)$ in the $2$-dimensional
subspace of $H^1(\cO_X)\oplus H^0(K_X)$ spanned by $(b,e)$,
$(c,f)$. We have only to check that the isomorphism classes of the
Higgs bundles given by (\ref{transition for tcE 1}), (\ref{Higgs
field of tcE 1}), (\ref{transition for tcE 2}) and (\ref{Higgs
field of tcE 2}) depend not on the particular choice of
$(b,e)$, $(c,f)$ but only on the points in
$\pp^{2g-4}$-bundle over $\Gr^\omega(2,2g)$. By Lemma
\ref{indep2}, the isomorphism classes of Higgs bundles given by
(\ref{transition for tcE 2}) and (\ref{Higgs field of tcE 2}) are
independent of the choice of $(b,e)$, $(c,f)$. On the other hand, it follows from Lemma \ref{indep1}-$(b)$ that the
isomorphism classes of the Higgs bundles given by (\ref{transition
for tcE 1}) and (\ref{Higgs field of tcE 1}) are independent of
the choice of $(g,h)$, $(b,e)$, $(c,f)$. Thus, there exists a morphism $\rho_{\s}:\bK_{\s}\to\bS$ such
that $\rho_{\s}\circ f_{\s}=\rho$.

We next show that $\rho_{\s}$ factors through $\bK_{\e}$. By Proposition \ref{contraction}, the
morphism $f_{\e}:\bK_{\s}\to\bK_{\e}$ is the contraction of the
fiber $\pp^5$ of a $\pp^5$-bundle over $\Gr^\omega(3,2g)$
and general points of a fiber give rise to a rank $4$ Higgs bundle given by
(\ref{transition for cE_1|_{z=0}}) and (\ref{Higgs field of
cE_1|_{z=0}}). By Lemma \ref{indep1}-$(a)$, the isomorphism
classes of the Higgs bundles given by (\ref{transition for
cE_1|_{z=0}}) and (\ref{Higgs field of cE_1|_{z=0}}) depend only
on the $3$-dimensional subspace spanned by $(a,d),(b,e),(c,f)$.
Thus $\rho_{\s}$ is constant along the fibers of $f_{\e}$. Applying Riemann's extension theorem again, we obtain a morphism
$\rho_{\e}:\bK_{\e}\to\bS$ such that $\rho_{\e}\circ f=\rho$.
\end{proof}

\subsection{A proof of Theorem \ref{main}}
We precisely refer to the original version of Zariski's main theorem.

\begin{prop}[Main Theorem of \cite{Zar43}]\label{ZMT}
If $W$ is an irreducible fundamental variety on $V$ of a birational correspondence $T$ between $V$ and $V'$ and if $T$ has no fundamental elements on $V'$, then-under the assumption that $V$ is locally normal at $W$-each irreducible component of the transform $T[W]$ is of higher dimension than $W$.
\end{prop}

By Proposition \ref{geometry of canonical subschemes} and Corollary \ref{geometry1 of 1-2 canonical subscheme},
$\rho(\tD^{(2)}-\tD^{(1)}\cup\tD^{(3)})=\bS_{1}-\bS_{2}$ is a
smooth divisor of $\bS-\rho(\tD^{(1)}\cup\tD^{(3)})=\bS-\bS_{2}$
lying over $T^*J/\zz_{2}-\zz_2^{2g}$. Then we have a morphism $g$
from $\bS-\rho(\tD^{(1)}\cup\tD^{(3)})$ to the blow-up of
$\bM_{2}-\zz_2^{2g}$ along $T^*J/\zz_{2}-\zz_2^{2g}$ which is
isomorphic to
$\bK-\tD^{(1)}\cup\tD^{(3)}=\bK_{\e}-f(\tD^{(1)}\cup\tD^{(3)})$. Thus, $\rho_{\e}$ mentioned in the proof of Proposition \ref{factors through} and $g$ are isomorphisms in codimension one. Note that $f(\tD^{(1)}\cup\tD^{(3)})$ has codimension $3$ in $\bK_{\e}$.

Assume that $\bK_{\e}$ is not isomorphic to $\bS$. Since $\bS$ is smooth and $\bS_{2}$ is the fundamental variety on $\bS$ of $\rho_{\e}$, it follows from Proposition \ref{ZMT} that $\bS_{2}$ has codimension greater than $3$. Further, since $\bK_{\e}$ is smooth, $f(\tD^{(1)}\cup\tD^{(3)})$ is the fundamental variety on $\bK_{\e}$ of $g$, it follows from Proposition \ref{ZMT} that $\bS_{2}$ has codimension less than $3$. Hence to avoid a contradiction, $\bK_{\e}$ should be isomorphic to $\bS$.


\begin{thebibliography}{BCHM10}

\bibitem{Bal} V. Balaji. \textit{Cohomology of certain moduli spaces of vector bundles}, Proc. Indian Acad. Sci. Math. Sci. \textbf{98} (1988), no.1, 1--24.

\bibitem{BNR} A. Beauville, M.S. Narasimhan and S. Ramanan. \textit{Spectral curves and the generalised theta divisor}, J. Reine Angew. Math. \textbf{398} (1989), 169--179.

\bibitem{BPG} S.B. Bradlow, O.G. Prada and P.B. Gothen. \textit{Moduli spaces of holomorphic triples over compact
Riemann surfaces}. Math. Ann. \textbf{328} (2004), 299--351.

\bibitem{BR} I. Biswas and S. Ramanan. \textit{An infinitesimal study of the moduli of Hitchin pairs}. J. Lond. Math. Soc. (2) \textbf{49} (1994), 219--231.

\bibitem{GM} W.M. Goldman and J.J. Millson. \textit{The deformation theory of representations of fundamental groups of compact K\"{a}hler manifolds}

\bibitem{H} T. Hausel. \textit{Compactification of moduli of Higgs bundles}, J. Reine Angew. Math. \textbf{503} (1998), 169--192.

\bibitem{Hart10} R. Hartshorne, Deformation Theory, {Graduate Texts in Mathematics} {\bf257}, Springer-Verlag, New York (2010).

\bibitem{HP} T. Hausel and C. Pauly. \textit{Prym varieties of spectral covers}, Geom. Topol. \textbf{16} (2012), no.3, 1609--1638.

\bibitem{HT} T. Hausel and M. Thaddeus. \textit{Generators for the cohomology ring of the moduli space of rank 2 Higgs bundles}, Proc. Lond. Math. Soc. (3) \textbf{88} (2004), no.3, 632--658.

\bibitem{Hit} N.J. Hitchin. \textit{The self-duality equations on a Riemann surface}, Proc.
Lond. Math. Soc. (3) \textbf{55} (1987), 59--126.

\bibitem{HL} D. Huybrechts and M. Lehn.
\textit{The Geometry of moduli spaces of sheaves.} A Publication
of the Max-Planck-Institut f\"{u}r Mathematik, Bonn, 1997.

\bibitem{KL} Y.-H. Kiem and J. Li. \textit{Desingularizations of the moduli space of rank 2 bundles over a curve}, Math. Ann. \textbf{330} (2004), 491--518.

\bibitem{K2} F. Kirwan. \textit{Partial desingularisations of quotients of nonsingular
varieties and their {B}etti numbers.} Ann. of Math. \textbf{122} (1985), 41--85.

\bibitem{KY} Y.-H. Kiem and S.-B. Yoo. \textit{The Stringy E-function of the moduli space of Higgs bundles with trivial determinant}, Math. Nachr. \textbf{281} (2008), no.6, 817--838.

\bibitem{Me} M. Mehta. \textit{Birational equivalence of Higgs moduli}. Internat. J. Math. \textbf{16} (2005), no.4, 365--386.

\bibitem{Mum76} D. Mumford. Algebraic Geometry. I. Complex Projective Varieties, Springer-Verlag (1976).

\bibitem{Nit} N. Nitsure. \textit{Moduli space of semistable pairs on a curve}, Proc. Lond. Math. Soc. (3) \textbf{62} (1991), 275--300.

\bibitem{OGr} K.G. O'Grady. \textit{Desingularized moduli spaces of sheaves on a K3}, J. Reine Angew. Math. \textbf{512} (1999), 49--117.

\bibitem{OGr97} K.G. O'Grady. \textit{Desingularized moduli spaces of sheaves on a K3, I}, arXiv:alg-geom/9708009, https://arxiv.org/pdf/alg-geom/9708009.pdf.

\bibitem{Se77} C.S. Seshadri. \textit{Desingularisation of moduli varieties of vector bundles on curves}, Proc. Int. Symp. on Algebraic Geometry, Kyoto, 155--184, 1977.

\bibitem{Se82} C.S. Seshadri. \textit{Fibr\'{e}s vectoriels sur les courbes
alg\'{e}briques}, June 1980. Asterisque, \textbf{96},
Soci\'{e}t\'{e} Math\'{e}matique de France, 1982.

\bibitem{Sim} C.T. Simpson. \textit{Moduli of representations of the fundamental group of a smooth projective
variety I and II}, Inst. Hautes Etudes Sci. Publ. Math. No. \textbf{79} (1994), 47--129 and \textbf{80} (1994), 5--79.

\bibitem{Yok1} K. Yokogawa. \textit{Compactification of moduli of parabolic sheaves and moduli of parabolic Higgs
sheaves}, J. Math. Kyoto Univ. \textbf{33-2} (1993), 451--504.

\bibitem{Yok2} K. Yokogawa. \textit{Infinitesimal deformation of parabolic Higgs
sheaves}, Internat. J. Math. \textbf{6} (1995), No. 1, 125--148.

\bibitem{Zar43} O. Zariski. \textit{Foundations of a general theory of birational correspondences}, Trans. Amer. Math. Soc. \textbf{53} (1943), 490–-542.

\end{thebibliography}
\end{document}